\documentclass[journal]{IEEEtran}
\usepackage[utf8]{inputenc}
\usepackage{amsmath}
\usepackage{amsfonts}
\usepackage{amssymb}
\usepackage{bbm}
\usepackage{amsthm}
\usepackage{mathtools}
\usepackage{float}
\usepackage{tkz-graph}
\usepackage{enumitem}
\usepackage{setspace}
\usepackage{graphicx}
\usepackage{color}
\usepackage{pgfplots}
\usepackage{cite}
\usepackage{subdepth}
\usepackage{verbatim} 
\usepackage{tikz}
\usetikzlibrary{arrows,positioning,decorations.pathreplacing,shapes,arrows,positioning,backgrounds}
\usepackage{algpseudocode,algorithm,algorithmicx}

\algrenewcommand\algorithmicrequire{\textbf{Input:}}
\algrenewcommand\algorithmicensure{\textbf{Output:}}
\DeclarePairedDelimiter\abs{\lvert}{\rvert}%
\newtheorem{theorem}{Theorem}
\newtheorem{lemma}{Lemma}
\newtheorem{problem}{Problem}
\newtheorem{proposition}{Proposition}
\newtheorem{corollary}{Corollary}
\newtheoremstyle{example}
{3pt}
{3pt}
{}
{}
{\upshape \bfseries}
{.}
{.5em}
{}
\theoremstyle{example}
\newtheorem{example}{Example}
\newtheorem{remark}{Remark}
\newtheorem{definition}{Definition}
\newenvironment{myproof}[1][\proofname]{\proof[\bf #1]}{\endproof}
\def\*#1{\mathbf{#1}}
\tikzstyle{bball} = [circle,shading=ball, ball color=black!100!white,
minimum size=0.4cm]
\tikzstyle{wball} = [circle,shading=ball, ball color=white!100!black,
minimum size=0.4cm]
\tikzset{VertexStyle1/.style = {shape = circle,
		color=black,
		fill=white!96!black,
		minimum size=0.5cm,
		text = black,
		inner sep = 2pt,
		outer sep = 1pt,
		minimum size = 0.55cm}
}		
\tikzset{VertexStyle2/.style = {shape = circle,
		color=black,
		fill=black!96!white,
		minimum size=0.5cm,
		text = white,
		inner sep = 2pt,
		outer sep = 1pt,
		minimum size = 0.55cm}
}
\DeclareMathOperator{\rank}{rank}

\newcommand{\adj}{\operatorname{adj}}
\newcommand \mc {\mathcal}

\newcommand \inv {^{-1}}

\begin{document} 
\title{Necessary and Sufficient Topological Conditions for Identifiability of Dynamical Networks}
\author{Henk J. van Waarde, Pietro Tesi, and M. Kanat Camlibel
	\thanks{Henk van Waarde and Kanat Camlibel are with the Bernoulli Institute for Mathematics, Computer Science and Artificial Intelligence, Faculty of Science and Engineering, University of Groningen, P.O. Box 407, 9700 AK Groningen, The Netherlands. Henk van Waarde is also with the Engineering and Technology Institute Groningen, Faculty of Science and Engineering, University of Groningen, 9747 AG Groningen, The Netherlands. Pietro Tesi is with the Engineering and Technology Institute Groningen and also with the Department of Information Engineering, University of Florence, 50139 Florence, Italy. Email: {\tt\footnotesize h.j.van.waarde@rug.nl, m.k.camlibel@rug.nl, p.tesi@rug.nl, pietro.tesi@unifi.it}. 
	}%
}

\maketitle

\begin{abstract}
This paper deals with dynamical networks for which the relations between node signals are described by proper transfer functions and external signals can influence each of the node signals. We are interested in graph-theoretic conditions for identifiability of such dynamical networks, where we assume that only a subset of nodes is measured but the underlying graph structure of the network is known. This problem has recently been investigated from a \emph{generic} viewpoint. Roughly speaking, generic identifiability means that the transfer functions in the network can be identified for ``\emph{almost all}" network matrices associated with the graph. In this paper, we investigate the stronger notion of identifiability \emph{for all} network matrices. To this end, we introduce a new graph-theoretic concept called the graph simplification process. Based on this process, we provide necessary and sufficient topological conditions for identifiability. Notably, we also show that these conditions can be verified by polynomial time algorithms. Finally, we explain how our results generalize existing sufficient conditions for identifiability.
\end{abstract}


\IEEEpeerreviewmaketitle

\section{Introduction}

\IEEEPARstart{N}{etworks} of dynamical systems appear in a variety of domains, including power systems, robotic networks, and aerospace systems \cite{Mesbahi2010}. In this paper, we consider a dynamical network model in which the relations between node signals are modelled by proper transfer functions. Such network models have received much attention in recent years, see e.g. \cite{vandenHof2013,Dankers2014,Hendrickx2018,Weerts2018,vanWaarde2018}. 

The interconnection structure of a dynamical network can be represented by a directed graph, where vertices (or nodes) represent scalar signals, and edges correspond to transfer functions connecting different node signals. We will assume that the underlying graph (i.e., the topology) of the dynamical network is \emph{known}. We remark that the related problem of \emph{topology identification} has also been studied, see e.g. \cite{Goncalves2008,Yuan2011,Nabi-Abdolyousefi2012,Materassi2012,Shahrampour2015,vanWaarde2017a}. 

We are interested in conditions for identifiability of dynamical networks. Identifiability is a fundamental property of a model set that guarantees that a unique (network) model can be identified, given informative data. Thus identifiability can be thought of as a prerequisite for identification: if identifiability does not hold then it is impossible to uniquely determine a network model, irrespective of the particular identification method and the experimental conditions. 

In the literature, several methods have been proposed for network identification \cite{vandenHof2013,Dankers2014,Weerts2018b,Haber2014}, these methods all exploit the structure of the network. For instance, a prediction error method was considered in \cite{Weerts2018b}, where consistency and minimum variance properties were proven under the assumption that the network is identifiable, the disturbances are filtered white noise, and the inputs are persistently exciting and uncorrelated with the disturbances. Another work \cite{Haber2014} considers subspace identification of networks with a path graph topology. As we will see, the structure of the network plays a fundamental role also with respect to the question of identifiability.

We follow the setup of \cite{Hendrickx2018}, where all network nodes can be externally excited, but only a subset of nodes can be measured. Within this setup, we are interested in two identifiability problems. Firstly, we want to find conditions under which the transfer functions from a given node to its out-neighbours are identifiable. Secondly, we wonder under which conditions the transfer functions of all edges in the network are identifiable. In particular, our aim is to find \emph{graph-theoretic} conditions for the above problems, that is, conditions in terms of the network structure and the locations of measured nodes. Such conditions based on the network topology are desirable since they give insight on the types of network structures that allow unique identification, and in addition may aid in the \emph{selection} of measured nodes. Graph-theoretic methods have also been succesfully applied to assess other system-theoretic properties like structural controllability \cite{Liu2011,Chapman2014,Monshizadeh2014,vanWaarde2017b,Jia2019} and fault detection \cite{Rapisarda2015,deRoo2015}.

Identifiability of dynamical networks is an active research area, see e.g. \cite{Hendrickx2018,Weerts2018,vanWaarde2018,Weerts2016,Adebayo2012,Hayden2017,Nabavi2016,vanWaarde2018b} and the references therein. The papers that are most closely related to the work presented here are \cite{Nabavi2016}, \cite{vanWaarde2018b}, \cite{Hendrickx2018}, and \cite{vanWaarde2018}, in which identifiability is also considered from \emph{graph-theoretic} perspective. In \cite{Nabavi2016} and \cite{vanWaarde2018b}, sufficient graph-theoretic conditions for identifiability have been presented for a class of \emph{state-space} systems. 

In \cite{Hendrickx2018}, graph-theoretic conditions have been established for \emph{generic} identifiability. That is, conditions were given under which transfer functions in the network can be identified for \emph{``almost all"} network matrices associated with the graph. The authors of \cite{Hendrickx2018} show that generic identifiability is equivalent to the existence of certain \emph{vertex-disjoint paths}, which yields ele\-gant conditions for generic identifiability. 

Inspired by the work in \cite{Hendrickx2018}, we are interested in graph-theoretic conditions for a stronger notion, namely identifiability \emph{for all} network matrices associated with the graph, a notion often referred to as \emph{global identifiability}. This problem is motivated by the fact that, although generic identifiability guarantees identifiability for almost all network matrices, there are meaningful examples of network matrices that are not contained in this set of almost all network matrices. As a consequence, a situation may arise in which the system under consideration is not identifiable, even though the conditions for generic identifiability are satisfied. For an example of such a situation, we refer to Section \ref{sectionproblem}. On the other hand, if the conditions derived in this paper are satisfied, then it is guaranteed that the network is identifiable \emph{for all} network matrices associated with the graph. 

The contributions of this paper are the following.
\begin{enumerate}
	\item We introduce the so-called \emph{graph simplification process}. Based on this process, we provide necessary and sufficient conditions for the left-invertibility of certain network-related transfer matrices. 
	\item Using the fact that identifiability is characterized by the left-invertibility of transfer matrices \cite{Hendrickx2018}, \cite{vanWaarde2018}, we provide necessary and sufficient graph-theoretic conditions for identifiability based on graph simplification. We also show that these conditions can be verified by polynomial time algorithms.
	\item We compare our results with the sufficient topological conditions for identifiability based on constrained vertex-disjoint paths \cite{vanWaarde2018}. In particular, we show that the results obtained in this paper generalize those in \cite{vanWaarde2018}.
\end{enumerate}

This paper is organized as follows. In Section \ref{sectionpreliminaries} we discuss the preliminaries that are used throughout this paper. Subsequently, in Section \ref{sectionproblem} we state and motivate the problem. Next, in Section \ref{sectionrank} we recall rank conditions for identifiability. Sections \ref{sectiongraphsimplification} and \ref{sectionidentifiabilitysimplification} contain our main results. In Section \ref{sectiongraphsimplification} we introduce the graph simplification process and show its relation to the left-invertibility of transfer matrices. Subsequently, in Section \ref{sectionidentifiabilitysimplification} we provide graph-theoretic conditions for identifiability. Our main results are compared to previous work in Section \ref{sectionresultscvdp}. Finally, Section \ref{sectionconclusions} contains our conclusions.

\section{Preliminaries}
\label{sectionpreliminaries}

We denote the set of natural numbers by $\mathbb{N}$, real numbers by $\mathbb{R}$, and complex numbers by $\mathbb{C}$. The set of real $m \times n$ matrices is denoted by $\mathbb{R}^{m \times n}$. The $n \times n$ identity matrix is denoted by $I_n$. If its dimension is clear, we simply write $I$.  

\subsection{Rational functions and rational matrices}
Consider a scalar variable $z$ and a rational function $f(z) = \frac{p(z)}{q(z)}$, where $p(z)$ and $q(z)$ are real polynomials and $q$ is nonzero. The function $f$ is \emph{proper} if the degree of $p(z)$ is less than or equal to the degree of $q(z)$. We say $f$ is \emph{strictly proper} if the degree of $p(z)$ is less than the degree of $q(z)$. An $m \times n$ matrix $A(z)$ is called \emph{rational} if its entries are rational functions in the variable $z$. In addition, $A(z)$ is \emph{proper} if its entries are proper rational functions. We omit the argument $z$ whenever the dependency of $A$ on $z$ is clear from the context. The \emph{normal rank} of $A(z)$ is defined as $\max_{\lambda \in \mathbb{C}} \rank A(\lambda)$ and denoted by $\rank A(z)$, with slight abuse of notation. We say $A(z)$ is \emph{left-invertible} if $\rank A(z) = n$. We denote the $(i,j)$-th entry of a matrix $A$ by $A_{ij}$. Moreover, the $j$-th column of $A$ is given by $A_{\bullet j}$. More generally, let $\mathcal{M} \subseteq \{1,2,\dots,m\}$ and $\mathcal{N} \subseteq \{1,2,\dots,n\}$. Then, $A_{\mathcal{M},\mathcal{N}}$ denotes the submatrix of $A$ containing the rows of $A$ indexed by $\mathcal{M}$ and the columns of $A$ indexed by $\mathcal{N}$. Next, consider the case that $A$ is square, i.e., $m = n$. The \emph{determinant} of $A$ is denoted by $\det A$, while the \emph{adjugate} of $A$ is denoted by $\adj A$. A \emph{principal submatrix} of $A$ is a submatrix $A_{\mathcal{M},\mathcal{M}}$, where $\mathcal{M} \subseteq \{1,2,\dots,m\}$. The determinant of $A_{\mathcal{M},\mathcal{M}}$ is called a \emph{principal minor} of $A$. 

\subsection{Graph theory}
\label{sectiongraphtheory}
Let $\mathcal{G} = (\mathcal{V},\mathcal{E})$ be a directed graph, with vertex (or node) set $\mathcal{V} = \{1,2,\dots,n\}$ and edge set $\mathcal{E} \subseteq \mathcal{V} \times \mathcal{V}$. The graphs considered in this paper are \emph{simple}, i.e., without self-loops and with at most one edge from one node to another. Consider an edge $(i,j) \in \mathcal{E}$. Then $(i,j)$ is called an \emph{outgoing} edge of node $i \in \mathcal{V}$ and $j$ is called an \emph{out-neighbour} of $i \in \mathcal{V}$. The set of out-neighbours of $i$ is denoted by $\mathcal{N}^+_i$. Similarly, $(i,j)$ is called an \emph{incoming} edge of $j \in \mathcal{V}$ and node $i$ is called an \emph{in-neighbour} of $j$. The set of in-neighbours of node $j$ is denoted by $\mathcal{N}_j^-$. For any subset $\mathcal{S} = \{v_1,v_2,\dots,v_s\} \subseteq \mathcal{V}$ we define the $s \times n$ matrix $P(\mathcal{V};\mathcal{S})$ as $P_{ij} := 1$ if $j = v_i$, and $P_{ij} := 0$ otherwise. The complement of $\mathcal{S}$ in $\mathcal{V}$ is defined as $\mathcal{S}^c := \mathcal{V} \setminus \mathcal{S}$. Moreover, the cardinality of $\mathcal{S}$ is denoted by $|\mathcal{S}|$. A \emph{path} $\mathcal{P}$ is a set of edges in $\mathcal{G}$ of the form $\mathcal{P} = \{(v_i,v_{i+1}) \mid i = 1,2,\dots,k \} \subseteq \mathcal{E}$, where the vertices $v_1,v_2,\dots,v_{k+1}$ are \emph{distinct}. The vertex $v_1$ is called a \emph{starting node} of $\mathcal{P}$, while $v_{k+1}$ is the \emph{end node}. The cardinality of $\mathcal{P}$ is called the \emph{length} of the path.  A collection of paths $\mathcal{P}_1,\mathcal{P}_2,\dots,\mathcal{P}_l$ is called \emph{vertex-disjoint} if the paths have no vertex in common, that is, if for all distinct $i,j \in \{1,2,\dots,l\}$, we have that
\begin{equation*}
(u_i,w_i) \in \mathcal{P}_i, (u_j,w_j) \in \mathcal{P}_j \implies u_i,w_i,u_j,w_j \text{ are distinct}.
\end{equation*}
Let $\mathcal{U}, \mathcal{W} \subseteq \mathcal{V}$ be disjoint. We say there exists a path \emph{from} $\mathcal{U}$ \emph{to} $\mathcal{W}$ if there exist vertices $u \in \mathcal{U}$ and $w \in \mathcal{W}$ such that there exists a path in $\mathcal{G}$ with starting node $u$ and end node $w$. Similarly, we say there are $m$ vertex-disjoint paths from $\mathcal{U}$ to $\mathcal{W}$ if there exist $m$ vertex-disjoint paths\footnote{Such sets of vertex-disjoint paths have been studied in detail in \cite{Murota1987}, where they were called linkings.} in $\mathcal{G}$ with starting nodes in $\mathcal{U}$ and end nodes in $\mathcal{W}$. In the case that $\mathcal{U} \cap \mathcal{W} \neq \emptyset$, we say there exist $m$ vertex-disjoint paths from $\mathcal{U}$ to $\mathcal{W}$ if there are $\max\{0,m - \abs{\mathcal{U} \cap \mathcal{W}}\}$ vertex-disjoint paths from $\mathcal{U} \setminus \mathcal{W}$ to $\mathcal{W} \setminus \mathcal{U}$. Roughly speaking, this means that we count paths of ``length zero" from every node in $\mathcal{U} \cap \mathcal{W}$ to itself.

\section{Problem statement and motivation}
\label{sectionproblem}
Let $\mathcal{G} = (\mathcal{V},\mathcal{E})$ be a simple directed graph with vertex set $\mathcal{V} = \{1,2,\dots,n\}$ and edge set $\mathcal{E} \subseteq \mathcal{V} \times \mathcal{V}$. We associate with each node $i \in \mathcal{V}$ a scalar \emph{node signal} $w_i(t)$, an external \emph{excitation signal} $r_i(t)$ and a \emph{disturbance signal} $v_i(t)$. Then, we consider the following discrete-time dynamics:
\begin{equation*}
w_i(t) = \sum_{j \in \mathcal{N}_i^-} G_{ij}(q) w_j(t) + r_i(t) + v_i(t),
\end{equation*}
where $G_{ij}(z)$ is a scalar transfer function and $q$ denotes the forward shift operator defined by $q w_i(t) = w_i(t+1)$. By concatenation of the node signals, excitation signals and disturbance signals, we can write the dynamics of all nodes compactly as $w(t) = G(q) w(t) + r(t) + v(t)$, where $w, r$, and $v$ are $n$-dimensional vectors and $G(z)$ is a $n \times n$ rational matrix. 
In addition, we consider a measured output vector $y(t)$ of dimension $p$ that consists of the node signals of a subset $\mathcal{C} \subseteq \mathcal{V}$ of so-called \emph{measured nodes}. By defining an associated binary matrix $C$ as $C := P(\mathcal{V},\mathcal{C})$, we can write this output as $y(t) = C w(t)$.
Finally, by combining the equations for $w$ and $y$, we obtain the networked system
\begin{equation}
\label{system}
\begin{aligned}
w(t) &= G(q) w(t) + r(t) + v(t) \\
y(t) &= C w(t).
\end{aligned} 
\end{equation}
We call the matrix $G(z)$ the \emph{network matrix} and assume that it satisfies the following properties:
\begin{enumerate}[label=\textbf{P\arabic*.}]
	\item For all $i,j \in \mathcal{V}$, the entry $G_{ij}(z)$ is a proper rational (transfer) function.
	\item The function $G_{ij}(z)$ is nonzero if and only if $(j,i) \in \mathcal{E}$. A matrix $G(z)$ that satisfies this property is said to be \emph{consistent} with the graph $\mathcal{G}$.
	\item Every principal minor of $\lim_{z \to \infty} (I-G(z))$ is nonzero. This implies that the network model \eqref{system} is \emph{well-posed} in the sense of Definition 2.11 of \cite{Dankers2014}.
\end{enumerate}
A network matrix $G(z)$ satisfying Properties P1, P2, and P3 is called \emph{admissible}. The set of all admissible network matrices is denoted by $\mathcal{A}(\mathcal{G})$. 
\begin{remark}
	A continuous-time counterpart of \eqref{system} can be obtained by replacing $q$ by the differential operator, hence our results will also be applicable to continuous-time systems. Besides the model \eqref{system}, also \emph{state-space} network models have received much attention (see, e.g., \cite{Goncalves2008,Yuan2011,Nabi-Abdolyousefi2012,Hayden2017,vanWaarde2017a}). A state-space model with scalar node dynamics can be obtained from \eqref{system} by choosing the nonzero entries of $G$ as first-order functions \cite{Kivits2018}. Also more general state-space models can be found in the literature, where the node dynamics are described by general linear systems, see e.g \cite{Fuhrmann2015}. The model \eqref{system} cannot capture these dynamics since the node signals $w_i$ are assumed to be scalar. The extension to the non-scalar case is therefore of interest, and will be considered for future work.
\end{remark}
For the development of this paper, it is important to distinguish between the following two concepts:
	\begin{itemize}
		\item \emph{Identifiability}: this is a fundamental property of the set of models of the form \eqref{system} that captures under what conditions identification is conceptually possible. If this property is not satisfied, one cannot uniquely identify the dynamics, no matter which identification algorithm is used. Identifiability does not involve any use of data.    
		\item \emph{Identification}: this involves the development of numerical algorithms for identifying the system dynamics from data. If identifiability holds then identification can be successfully performed in different ways under hypotheses on the noise and the informativity of the data \cite{Ljung1999}. 
	\end{itemize} 
This paper focuses on characterizations of \emph{identifiability}. To explain what identifiability means in a network context, we first write \eqref{system} in input/output form as
\begin{equation*}
y(t) = C(I-G(q))\inv r(t) + \bar{v}(t),
\end{equation*}
where $\bar{v}(t) := C(I-G(q))\inv v(t)$. It is well-known that the transfer matrix $C(I-G(z))\inv$ from $r$ to $y$ can be obtained from $\{r(t),y(t)\}$-data, under suitable assumptions on $r$ and $v$ \cite{Ljung1999}. The question of network identifiability is then the following: which entries of $G(z)$ can be uniquely reconstructed from $C(I-G(z))\inv$? In this paper we restrict our attention to the identifiability of the transfer functions outgoing a given node $i$ (i.e., identifiability of a column of $G(z)$), and to the identifiability of the entire matrix $G(z)$. A standing assumption in our work is that we \emph{know} the graph structure $\mathcal{G}$ underlying the dynamical network.

In recent work \cite{Hendrickx2018}, \cite{Bazanella2017} the problem of identifiability has been considered from \emph{generic} viewpoint. Graph-theoretic conditions were given under which certain entries of $G(z)$ can be uniquely reconstructed from $C(I-G(z))\inv$ \emph{for almost all} network matrices $G$ consistent with the graph. For a formal definition of generic identifiability we refer to Definition 1 of \cite{Hendrickx2018}. Here, we will informally illustrate the approach of \cite{Hendrickx2018}. We will use the shorthand notation $T(z;G) := (I - G(z))\inv$. This means that the transfer matrix from $r$ to $y$ equals $CT$. 
\begin{example}
	\label{example1}
	Consider the graph $\mathcal{G} = (\mathcal{V},\mathcal{E})$ in Figure \ref{fig:graph1}. We assume that the node signals of nodes $4$ and $5$ can be measured, that is, $\mathcal{C} = \{4,5\}$. Suppose that we want to identify the transfer functions from node $1$ to its out-neighbours, i.e., the transfer functions $G_{21}(z)$ and $G_{31}(z)$. According to Corollary 5.1 of \cite{Hendrickx2018}, this is possible if and only if there exist two vertex-disjoint paths from $\mathcal{N}^+_1$ to $\mathcal{C}$. Note that this is the case in this example, since the edges $(2,4)$ and $(3,5)$ are two vertex-disjoint paths.	To see why we can generically identify the transfer functions $G_{21}$ and $G_{31}$, we compute $CT$ as:
	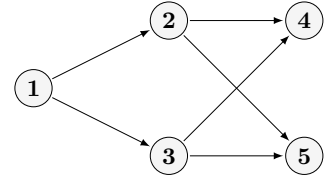
\begin{figure}[h!]
		\centering
		\scalebox{0.9}{
			\begin{tikzpicture}[scale=1]
			\node[VertexStyle1,draw] (1) at (0,0) {$\bf 1$};
			\node[VertexStyle1,draw] (2) at (2,1) {$\bf 2$};
			\node[VertexStyle1,draw] (3) at (2,-1) {$\bf 3$};
			\node[VertexStyle1,draw] (4) at (4,1) {$\bf 4$};
			\node[VertexStyle1,draw] (5) at (4,-1) {$\bf 5$};
			\draw[-latex] (1) -- (2);
			\draw[-latex] (1) -- (3);
			\draw[-latex] (2) -- (4);
			\draw[-latex] (2) -- (5);
			\draw[-latex] (3) -- (4);
			\draw[-latex] (3) -- (5);
			\end{tikzpicture}
		}
		\caption{Graph used in Example \ref{example1}.}
		\label{fig:graph1}
	\end{figure}
	\begin{equation*}
	CT = \begin{pmatrix}
	G_{42} G_{21} + G_{43} G_{31} & G_{42} & G_{43} & 1 & 0 \\
	G_{52} G_{21} + G_{53} G_{31} & G_{52} & G_{53} & 0 & 1
	\end{pmatrix},
	\end{equation*}
	where we omit the argument $z$. Clearly, we can uniquely obtain the transfer functions $G_{42}, G_{43}, G_{52}$, and $G_{53}$ from $CT$. Moreover, the transfer matrices $G_{21}$ and $G_{31}$ satisfy 
	\begin{equation}
	\label{equationex1}
	\begin{pmatrix}
	G_{42} & G_{43} \\
	G_{52} & G_{53}
	\end{pmatrix}
	\begin{pmatrix}
	G_{21} \\ G_{31}
	\end{pmatrix}
	= \begin{pmatrix}
	T_{41} \\ T_{51}
	\end{pmatrix}.
	\end{equation} 
	Equation \eqref{equationex1} has a unique solution in the unknowns $G_{21}$ and $G_{31}$ if $G_{42} G_{53} - G_{43} G_{52} \neq 0$, which means that we can identify $G_{21}$ and $G_{31}$ for ``almost all" $G$ consistent with $\mathcal{G}$. 
\end{example}

As mentioned before, the approach based on vertex-disjoint paths \cite{Hendrickx2018} gives necessary and sufficient conditions for \emph{generic} identifiability. This implies that for some network matrices $G$, it might be impossible to identify the transfer functions, even though the path-based conditions are satisfied. For instance, in Example \ref{example1} we cannot identify $G_{21}$ and $G_{31}$ if the network matrix $G$ is such that $G_{42} = G_{43} = G_{52} = G_{53}$. Nonetheless, scenarios in which some (or all) transfer functions in the network are equal occur frequently, for example in the study of undirected (electrical) networks \cite{dorfler2018}, in unweighted consensus networks \cite{Olfati-Saber2007}, and in the study of Cartesian products of graphs \cite{Chapman2013}. Therefore, instead of generic identifiability, we are interested in graph-theoretic conditions that guarantee identifiability \emph{for all} admissible network matrices. Such a problem might seem like a simple extension of the work on generic identifiability \cite{Hendrickx2018}. However, to analyze \emph{strong structural} network properties (for all network matrices), we typically need completely different graph-theoretic tools than the ones used in the analysis of \emph{generic} network properties. For instance, in the literature on \emph{controllability}, weak structural controllability is related to maximal matchings \cite{Liu2011}, while strong structural controllability is related to zero forcing sets \cite{Monshizadeh2014} and constrained matchings \cite{Chapman2013}. To make the problem of this paper more precise, we state a few definitions. First, we are interested in conditions under which all transfer functions from a node $i$ to its out-neighbours $\mathcal{N}^+_i$ are identifiable (for any admissible network matrix $G \in \mathcal{A}(\mathcal{G})$). If this is the case, we say $(i,\mathcal{N}^+_i)$ is \emph{globally identifiable}, or simply $(i,\mathcal{N}^+_i)$ is identifiable for short.

\begin{definition}
	\label{def1}
	Consider a directed graph $\mathcal{G} = (\mathcal{V},\mathcal{E})$ and let $i \in \mathcal{V}$ and $\mathcal{C} \subseteq \mathcal{V}$. Moreover, define $C = P(\mathcal{V},\mathcal{C})$. We say $(i,\mathcal{N}^+_i)$ is \emph{(globally) identifiable} from $\mathcal{C}$ if the implication
	\begin{equation*}
	CT(z;G) = CT(z;\bar{G}) \implies G_{\bullet i}(z) = \bar{G}_{\bullet i}(z)
	\end{equation*}
	holds for all $G(z), \bar{G}(z) \in \mathcal{A}(\mathcal{G})$.
\end{definition}

In addition, we are interested in conditions under which the \emph{entire} network matrix $G$ can be identified. If this is the case, we say the graph $\mathcal{G}$ is (globally) identifiable.

\begin{definition}
	Consider a directed graph $\mathcal{G} = (\mathcal{V},\mathcal{E})$ and let $\mathcal{C} \subseteq \mathcal{V}$ and $C = P(\mathcal{V},\mathcal{C})$. We say $\mathcal{G}$ is \emph{(globally) identifiable} from $\mathcal{C}$ if the implication 
	\begin{equation*}
	CT(z;G) = CT(z;\bar{G}) \implies G(z) = \bar{G}(z)
	\end{equation*}
	holds for all $G(z), \bar{G}(z) \in \mathcal{A}(\mathcal{G})$.
\end{definition}

The main goals of this paper are to find graph-theoretic conditions for identifiability of $(i,\mathcal{N}_i^+)$ and $\mathcal{G}$.

\begin{problem}
\label{problem1}
Consider a directed graph $\mathcal{G} = (\mathcal{V},\mathcal{E})$ with measured nodes $\mathcal{C} \subseteq \mathcal{V}$. Provide necessary and sufficient graph-theoretic conditions under which, respectively, $(i,\mathcal{N}_i^+)$ and $\mathcal{G}$ are identifiable from $\mathcal{C}$. 
\end{problem}

Graph-theoretic conditions for global identifiability are attractive for two reasons. First, such conditions will give insight on the types of graph structures that allow identification. Secondly, they allow us to select measured nodes guaranteeing identifiability \emph{before} collecting data. 
To deal with Problem \ref{problem1}, we make use of rank conditions for identifiability which we will recall in Section \ref{sectionrank}. To verify such rank conditions, we introduce a novel graph-theoretic concept called the \emph{graph simplification process} in Section \ref{sectiongraphsimplification}.
 
\section{Rank conditions for identifiability}
\label{sectionrank}
First, we review some of the conditions for identifiability in terms of the normal rank of transfer matrices. For the proofs of all results in this section, we refer to \cite{vanWaarde2018}. Recall from Section \ref{sectionpreliminaries} that $T_{\mathcal{C}, \mathcal{N}^+_i}(z;G)$ denotes the submatrix of $T$ formed by taking the rows of $T$ indexed by $\mathcal{C}$ and the columns of $T$ corresponding to $\mathcal{N}^+_i$. This means that $T_{\mathcal{C}, \mathcal{N}^+_i}(z;G)$ is a submatrix of the transfer matrix $CT(z;G)$ from $r$ to $y$, obtained by selecting the columns corresponding to $\mathcal{N}_i^+$. The following lemma (Lemma 5 of \cite{vanWaarde2018}) asserts that identifiability of $(i,\mathcal{N}^+_i)$ is equivalent to a rank condition on the matrix $T_{\mathcal{C}, \mathcal{N}^+_i}(z;G)$. 

\begin{lemma}
\label{lemma1}
Consider a directed graph $\mathcal{G} = (\mathcal{V},\mathcal{E})$, let $i \in \mathcal{V}$, and $\mathcal{C} \subseteq \mathcal{V}$. Then, $(i,\mathcal{N}^+_i)$ is identifiable from $\mathcal{C}$ if and only if $\rank T_{\mathcal{C}, \mathcal{N}^+_i}(z;G) = | \mathcal{N}^+_i |$ for all $G(z) \in \mathcal{A}(\mathcal{G})$.
\end{lemma}

As an immediate consequence of Lemma \ref{lemma1}, we find conditions for the identifiability of $\mathcal{G}$ based on the normal rank of transfer matrices. This is stated in the following corollary.

\begin{corollary}
\label{cor1}
	Consider a directed graph $\mathcal{G} = (\mathcal{V},\mathcal{E})$ and let $\mathcal{C} \subseteq \mathcal{V}$. Then, $\mathcal{G}$ is identifiable from $\mathcal{C}$ if and only if $\rank T_{\mathcal{C}, \mathcal{N}^+_i}(z;G) = | \mathcal{N}^+_i |$ for all $i \in \mathcal{V}$ and all $G(z) \in \mathcal{A}(\mathcal{G})$.
\end{corollary}

Although Lemma \ref{lemma1} and Corollary \ref{cor1} give necessary and suffi\-cient conditions for the identifiability of respectively $(i,\mathcal{N}_i^+)$ and $\mathcal{G}$, these conditions are limited since there is no obvious method to \emph{check} left-invertibility of $T_{\mathcal{C}, \mathcal{N}^+_i}(z;G)$ for an infinite number of matrices $G$. Therefore, one of the main results of this paper will be graph-theoretic conditions for the left-invertibility of $T_{\mathcal{W}, \mathcal{U}}(z;G)$, where $\mathcal{U}, \mathcal{W} \subseteq \mathcal{V}$ are any two subsets of vertices. These conditions will be introduced in the next section. 

\section{The graph simplification process}
\label{sectiongraphsimplification}

In this section we provide necessary and sufficient conditions for left-invertibility of $T_{\mathcal{W},\mathcal{U}}(z;G)$ for all $G(z) \in \mathcal{A}(\mathcal{G})$, where $\mathcal{U}, \mathcal{W} \subseteq \mathcal{V}$. Loosely speaking, the idea is to simplify the graph $\mathcal{G}$ and nodes $\mathcal{W}$ in such a way that checking left-invertibility becomes easy. To give the reader some intuition for the approach, we start with the following basic lemma, which asserts that $T_{\mathcal{W},\mathcal{U}}(z;G)$ is left-invertible if $\mathcal{U} \subseteq \mathcal{W}$. 

\begin{lemma}
	\label{lemmanecsys}
	Consider a directed graph $\mathcal{G} = (\mathcal{V},\mathcal{E})$ and let $\mathcal{U}, \mathcal{W} \subseteq \mathcal{V}$. If $\mathcal{U} \subseteq \mathcal{W}$ then $\rank T_{\mathcal{W},\mathcal{U}}(z;G) = \abs{\mathcal{U}}$ for all $G(z) \in \mathcal{A}(\mathcal{G})$.
\end{lemma}

The proof of Lemma \ref{lemmanecsys} is postponed to Appendix \ref{appendixlemmanecsys}. The condition $\mathcal{U} \subseteq \mathcal{W}$ considered in Lemma \ref{lemmanecsys} is clearly not necessary for left-invertibility. One can show this using the example $\mathcal{G} = (\mathcal{V},\mathcal{E})$, where $\mathcal{V} = \{1,2\}$, $\mathcal{E} = \{(1,2)\}$, and the subsets $\mathcal{U}$ and $\mathcal{W}$ are chosen as $\mathcal{U} = \{1\}$ and $\mathcal{W} = \{2\}$. However, the \emph{main idea} of the graph simplification process is to simplify $\mathcal{G}$ and to `move' the nodes in $\mathcal{W}$ closer to the nodes in $\mathcal{U}$ such that the condition $\mathcal{U} \subseteq \mathcal{W}$ possibly holds \emph{after} applying these operations. Of course, we cannot blindly modify the graph $\mathcal{G}$ since this would affect the left-invertibility of $T_{\mathcal{W},\mathcal{U}}(z;G)$. Instead, we will now state two lemmas in which we consider two different operations on $\mathcal{G}$ and $\mathcal{W}$ that \emph{preserve} left-invertibility of $T_{\mathcal{W},\mathcal{U}}(z;G)$. We emphasize that the graph operations are introduced for analysis purposes only. Indeed, since the condition of Lemma \ref{lemmanecsys} is simple to check, the graph operations should be seen as a \emph{tool} to check left-invertibility of the transfer matrix of a given \emph{fixed} graph $\mathcal{G}$. First, we state Lemma \ref{lemmaprocess1} which asserts that left-invertibility of $T_{\mathcal{W},\mathcal{U}}(z;G)$ is unaffected by the removal of the outgoing edges of $\mathcal{W}$.

\begin{lemma}
	\label{lemmaprocess1}
	Consider a directed graph $\mathcal{G} = (\mathcal{V},\mathcal{E})$ and let $\mathcal{U}, \mathcal{W} \subseteq \mathcal{V}$. Moreover, let $\bar{\mathcal{G}} = (\mathcal{V},\bar{\mathcal{E}})$ be the graph obtained from $\mathcal{G}$ by removing all outgoing edges of the nodes in $\mathcal{W}$. Then $\rank T_{\mathcal{W},\mathcal{U}}(z;G) = \abs{\mathcal{U}}$ for all $G(z) \in \mathcal{A}(\mathcal{G})$ if and only if $\rank T_{\mathcal{W},\mathcal{U}}(z;\bar{G}) = \abs{\mathcal{U}}$ for all $\bar{G}(z) \in \mathcal{A}(\bar{\mathcal{G}})$.
\end{lemma}

\begin{myproof}
	Let $G(z) \in \mathcal{A}(\mathcal{G})$. Relabel the nodes in $\mathcal{V}$ such that 
	\begin{equation}
	\label{matrixG}
	G = \begin{pmatrix}
	G_{\mathcal{R},\mathcal{R}} & G_{\mathcal{R},\mathcal{W}} \\
	G_{\mathcal{W},\mathcal{R}} & G_{\mathcal{W},\mathcal{W}}
	\end{pmatrix},
	\end{equation}
	where $\mathcal{R} := \mathcal{V} \setminus \mathcal{W}$ and the argument $z$ has been omitted. Define the matrix $\bar{G}$ as
	\begin{equation}
	\label{matrixGbar}
	\bar{G} = \begin{pmatrix}
	G_{\mathcal{R},\mathcal{R}} & 0 \\
	G_{\mathcal{W},\mathcal{R}} & 0
	\end{pmatrix}.
	\end{equation}
	The matrix $\bar{G}$ is an admissible matrix consistent with $\bar{\mathcal{G}}$, i.e., $\bar{G} \in \mathcal{A}(\bar{\mathcal{G}})$. To see this, note that $\bar{G}$ satisfies Property P1. Moreover, since all outgoing edges of nodes in $\mathcal{W}$ are removed in the graph $\bar{\mathcal{G}}$, the matrix $\bar{G}$ is consistent with $\bar{\mathcal{G}}$. Hence, $\bar{G}$ satisfies property P2. Finally, to see that $\bar{G}$ satisfies Property P3, note that any principal minor of 
	\begin{equation}
	\label{eqlim}
	\lim_{z \to \infty} \begin{pmatrix}
	I-G_{\mathcal{R},\mathcal{R}}(z) & 0 \\
	-G_{\mathcal{W},\mathcal{R}}(z) & I
	\end{pmatrix}
	\end{equation}
	is either 1 or equal to a principal minor of $\lim_{z \to \infty} (I-G_{\mathcal{R},\mathcal{R}}(z))$, which is nonzero by the assumption that $G$ is admissible. We conclude that $\bar{G} \in \mathcal{A}(\bar{\mathcal{G}})$. Next, by Proposition 2.8.7 of \cite{Bernstein2011}, the inverse of $I - G$ can be written as
	\begin{equation*}
	T = (I - G)\inv = 
	\begin{pmatrix}
	\ast & \ast \\
	S(G) G_{\mathcal{W},\mathcal{R}} (I-G_{\mathcal{R},\mathcal{R}})\inv & S(G)
	\end{pmatrix},
	\end{equation*}
	where $S(G) := (I-G_{\mathcal{W},\mathcal{W}} - G_{\mathcal{W},\mathcal{R}}(I - G_{\mathcal{R},\mathcal{R}})\inv G_{\mathcal{R},\mathcal{W}}))\inv$ denotes the inverse Schur complement of $I - G$. Using the same formula to compute the inverse of $I - \bar{G}$, we find 
	\begin{equation*}
	\bar{T} := (I - \bar{G})\inv = 
	\begin{pmatrix}
	\ast & \ast \\
	G_{\mathcal{W},\mathcal{R}} (I-G_{\mathcal{R},\mathcal{R}})\inv & I
	\end{pmatrix}.
	\end{equation*}
	The above expressions for $T$ and $\bar{T}$ imply that
	\begin{equation*}
	T_{\mathcal{W},\mathcal{U}} = S(G) \bar{T}_{\mathcal{W},\mathcal{U}},
	\end{equation*}
	and because $S(G)$ has full normal rank, we obtain
	\begin{equation}
	\label{rankequal}
	\rank T_{\mathcal{W},\mathcal{U}} = \rank \bar{T}_{\mathcal{W},\mathcal{U}}.
	\end{equation}
	Next, we use \eqref{rankequal} to prove the lemma. First, to prove the `if' statement, suppose that $\rank T_{\mathcal{W},\mathcal{U}}(z;\bar{G}) = \abs{\mathcal{U}}$ for all matrices $\bar{G} \in \mathcal{A}(\bar{\mathcal{G}})$. Let $G \in \mathcal{A}(\mathcal{G})$. Using $G$, construct the matrix $\bar{G} \in \mathcal{A}(\bar{\mathcal{G}})$ in \eqref{matrixGbar}. By hypothesis, $\rank T_{\mathcal{W},\mathcal{U}}(z;\bar{G}) = \abs{\mathcal{U}}$ and therefore we conclude from \eqref{rankequal} that $\rank T_{\mathcal{W},\mathcal{U}}(z;G) = \abs{\mathcal{U}}$.
	
	Subsequently, to prove the `only if' statement, suppose that $\rank T_{\mathcal{W},\mathcal{U}}(z;G) = \abs{\mathcal{U}}$ for all $G(z) \in \mathcal{A}(\mathcal{G})$. Consider any matrix $\bar{G}(z) \in \mathcal{A}(\bar{\mathcal{G}})$ and note that $\bar{G}$ can be written in the form \eqref{matrixGbar}. Next, we choose the matrices $G_{\mathcal{R},\mathcal{W}}$ and $G_{\mathcal{W},\mathcal{W}}$ such that the matrix $G$ in \eqref{matrixG} is consistent with the graph $\mathcal{G}$, and such that the nonzero entries of $G_{\mathcal{R},\mathcal{W}}$ and $G_{\mathcal{W},\mathcal{W}}$ are \emph{strictly proper} rational functions. This means that $G$ readily satisfies Properties P1 and P2 (see Section \ref{sectionproblem}). In fact, $G$ also satisfies P3. Indeed, since $\lim_{z\to \infty} (I-G(z))$ is given by \eqref{eqlim}, it follows that every principal minor of $\lim_{z\to \infty} (I-G(z))$ is either 1 or equal to a principal minor of $\lim_{z \to \infty} (I-G_{\mathcal{R},\mathcal{R}})$, which is nonzero by the hypothesis that $\bar{G}(z) \in \mathcal{A}(\bar{\mathcal{G}})$. We conclude that $G$ satisfies Properties P1, P2, and P3, equivalently, $G \in \mathcal{A}(\mathcal{G})$.
	By hypothesis, $\rank T_{\mathcal{W},\mathcal{U}}(z;G) = \abs{\mathcal{U}}$ and consequently, by \eqref{rankequal} we conclude that $\rank T_{\mathcal{W},\mathcal{U}}(z;\bar{G}) = \abs{\mathcal{U}}$. This proves the lemma.
\end{myproof}

\begin{remark}
	\label{remarkUincoming}
	In similar fashion as in the proof of Lemma \ref{lemmaprocess1}, we can prove that all \emph{incoming} edges of nodes in $\mathcal{U}$ can be removed without affecting the left-invertibility of $T_{\mathcal{W},\mathcal{U}}(z;G)$. 
\end{remark}

Inspired by Lemma \ref{lemmaprocess1}, we wonder what type of operations we can further perform on the graph $\mathcal{G}$ and nodes $\mathcal{W}$ without affecting left-invertibility of $T_{\mathcal{W},\mathcal{U}}(z;G)$. In what follows we will show that under suitable conditions it is possible to `move' the nodes in $\mathcal{W}$ closer to the nodes in $\mathcal{U}$. Here the notion of \emph{reachability} in graphs will play an important role. For a subset $\mathcal{U} \subseteq \mathcal{V}$ and a node $j \in \mathcal{V} \setminus \mathcal{U}$, we say $j$ is \emph{reachable} from $\mathcal{U}$ if there exists at least one path from $\mathcal{U}$ to $j$. By convention, if $j \in \mathcal{U}$ then $j$ is reachable from $\mathcal{U}$. In the following lemma, we will show that the rank of $T_{\mathcal{W},\mathcal{U}}(z;G)$ is unaffected if we replace a node $k \in \mathcal{W} \setminus \mathcal{U}$ by $j$, provided that $j$ is the \emph{only} in-neighbour of $k$ that is reachable from $\mathcal{U}$.
\begin{lemma}
	\label{lemmaprocess2}
	Consider a directed graph $\mathcal{G} = (\mathcal{V},\mathcal{E})$ and let $\mathcal{U},\mathcal{W} \subseteq \mathcal{V}$. Suppose that $k \in \mathcal{W} \setminus \mathcal{U}$ has exactly one in-neighbour $j \in \mathcal{N}_k^-$ that is reachable from $\mathcal{U}$. Then for all $G(z) \in \mathcal{A}(\mathcal{G})$, we have 
	\begin{equation*}
	\rank T_{\mathcal{W},\mathcal{U}}(z;G) = \rank T_{\mathcal{\bar{W}},\mathcal{U}}(z;G),
	\end{equation*}
	where $\bar{\mathcal{W}} := (\mathcal{W} \setminus \{k\}) \cup \{j\}$.
\end{lemma}

\begin{remark}
	We emphasize that Lemma \ref{lemmaprocess2} does not require node $k$ to have exactly one in-neighbour. In general, node $k$ may have multiple in-neighbours, but if exactly one of such neighbours is reachable from $\mathcal{U}$, we can apply Lemma \ref{lemmaprocess2}. The intuition of Lemma \ref{lemmaprocess2} is as follows: under the assumptions, all information from the nodes in $\mathcal{U}$ enters node $k$ \emph{via} node $j$. Therefore, choosing node $k$ or node $j$ as a node in $\mathcal{W}$ does not make any difference. An interesting special case is obtained when \emph{both} nodes $j$ and $k$ are contained in $\mathcal{W}$. In this case, we obtain $\bar{\mathcal{W}} = \mathcal{W} \setminus \{k\}$, that is, node $k$ can be removed from $\mathcal{W}$ without affecting the rank of $T_{\mathcal{W},\mathcal{U}}(z;G)$.
\end{remark}

\begin{myproof}[Proof of Lemma \ref{lemmaprocess2}]
	By Lemma \ref{lemmaprocess1}, we can assume without loss of generality that the nodes in $\mathcal{W}$ have no outgoing edges. Let $G(z) \in \mathcal{A}(\mathcal{G})$. In what follows we omit the dependence of $G$ on $z$ and the dependence of $T(z;G)$ on both $z$ and $G$. Consider a vertex $v \in \mathcal{U}$. Note that
	\begin{subequations}
		\begin{align}
		(I - G) T &= I \\
		\sum_{l = 1}^n (I - G)_{k l} T_{l v} &= 0, \label{sumGT}
		\end{align}
	\end{subequations}
	where $n := |\mathcal{V}|$ and \eqref{sumGT} follows from the fact that $k \in \mathcal{W} \setminus \mathcal{U}$ and $v \in \mathcal{U}$ are \emph{distinct}. Equation \eqref{sumGT} implies that 
	\begin{equation}
	\label{eqTkv}
	T_{kv} = \sum_{l \in \mathcal{N}_k^{-}} G_{k l} T_{l v}. 
	\end{equation}
	Note that $j \in \mathcal{N}_k^{-}$, but possibly $\mathcal{N}_k^-$ contains other vertices. We will now prove that for all these other vertices, the corresponding transfer function $T_{lv}$ equals zero. That is, $T_{lv} = 0$ for all $l \in \mathcal{N}_k^- \setminus \{j\}$. To see this, we first observe that there does not exist a path in $\mathcal{G}$ from $v$ to $l \in \mathcal{N}_k^- \setminus \{j\}$. Indeed, suppose that there is a path $\mathcal{P}$ from $v$ to $l$. Then this path cannot contain the edge $(j,k)$, since node $k \in \mathcal{W} \setminus \mathcal{U}$ does not have any outgoing edges. This implies that there exists a path $\mathcal{P} \cup (l,k)$ from $v$ to $k$ via node $l$. This is a contradiction since by hypothesis $j$ is the only in-neighbour of $k$ that is reachable from $\mathcal{U}$. Therefore, we conclude that there does not exist a path from $v$ to $l$. By Lemma 3 of \cite{vandenHof2013} we conclude that $T_{lv} = 0$. This means that \eqref{eqTkv} can be simplified as
	\begin{equation*}
	T_{kv} = G_{k j} T_{j v}.
	\end{equation*}
	Since $v \in \mathcal{U}$ is arbitrary, it follows that 
	\begin{equation*}
	T_{k,\mathcal{U}} = G_{kj} T_{j,\mathcal{U}}.
	\end{equation*}
	As $G_{kj} \neq 0$, we conclude that
	\begin{equation*}
	\rank T_{\mathcal{W},\mathcal{U}} = \rank T_{\bar{\mathcal{W}},\mathcal{U}},
	\end{equation*}
	where $\bar{\mathcal{W}} := (\mathcal{W} \setminus \{k\}) \cup \{j\}$. This proves the lemma.
\end{myproof}

From Lemma \ref{lemmaprocess1} and Lemma \ref{lemmaprocess2}, we see that (i) we can always remove the outgoing edges of nodes in $\mathcal{W}$ and (ii) we can `move' nodes in $\mathcal{W}$ closer to $\mathcal{U}$ under suitable conditions. Of course, since both  operations do not affect left-invertibility of $T_{\mathcal{W},\mathcal{U}}$, we can also apply these operations multiple times consecutively. Therefore, we introduce the following process to simplify the graph $\mathcal{G}$ and move the nodes in $\mathcal{W}$. The idea of this process is to apply the above operations to the graph until no more changes are possible. \\[2mm]
\noindent
\fbox{
	\parbox{0.465\textwidth}{
		\textbf{Graph simplification process:} \\ 
		\noindent
		Let $\mathcal{G} = (\mathcal{V},\mathcal{E})$ be a directed graph and let $\mathcal{U},\mathcal{W} \subseteq \mathcal{V}$. Consider the following two operations on the graph $\mathcal{G}$ and nodes $\mathcal{W}$.
		\begin{enumerate}
			\item Remove all outgoing edges of nodes in $\mathcal{W}$ from $\mathcal{G}$.
			\item If $k \in \mathcal{W} \setminus \mathcal{U}$ has exactly one in-neighbour $j \in \mathcal{N}_k^-$ that is reachable from $\mathcal{U}$, replace $k$ by $j$ in $\mathcal{W}$.
		\end{enumerate}
		Consecutively apply operations 1 and 2 on the graph $\mathcal{G}$ and nodes $\mathcal{W}$ until no more changes are possible.}} \\\\

\noindent
Clearly, the graph simplification process terminates after a finite number of applications of operations 1 and 2. Indeed, operation 1 can only be applied once in a row, and a node in $\mathcal{W} \setminus \mathcal{U}$ can be `moved' at most $\abs{\mathcal{V}}-1$ times which means that operation 2 can be applied only a finite number of times. In fact, it is attractive to apply the operations 1 and 2 in alternating fashion since the process will then terminate within $\abs{\mathcal{V}}$ operations of both types. This is due to the fact that if the outgoing edges of a node $j \in \mathcal{V}$ are removed, then we cannot apply operation 2 to replace a node $k$ by $j$. A graph obtained by applying the graph simplification process to $\mathcal{G}$ is called a \emph{derived graph}, which we denote by $\mathcal{D}(\mathcal{G})$. Similarly, we call a vertex set obtained by applying the graph simplification process to $\mathcal{W}$ a \emph{derived vertex set}, denoted by $\mathcal{D}(\mathcal{W})$. To stress the fact that $\mathcal{D}(\mathcal{G})$ and $\mathcal{D}(\mathcal{W})$ do not only depend on the graph $\mathcal{G}$ and set $\mathcal{W}$, but \emph{also} on the set $\mathcal{U}$, we say that $\mathcal{D}(\mathcal{G})$ and $\mathcal{D}(\mathcal{W})$ are a derived graph of $\mathcal{G}$ and derived vertex set of $\mathcal{W}$ \emph{with respect to} the set $\mathcal{U}$. We emphasize that derived graphs and derived vertex sets are not necessarily unique. In general, the derived graph and derived vertex set that are obtained from the graph simplification process depend on the \emph{order} in which the operations 1 and 2 are applied, and on the order in which operation 2 is applied to the nodes in $\mathcal{W}$. However, it turns out that the non-uniqueness of derived graphs and derived vertex sets is not a problem for the application (left-invertibility) we have in mind. In fact, we will show in Theorem \ref{maintheorem} that \emph{any} derived graph and derived vertex set will lead to the same conclusions about left-invertibility. 

\begin{remark}
	\label{remarkcomputation}
	In step 2 of the graph simplification process, we have to decide whether there exists a node $k \in \mathcal{W} \setminus \mathcal{U}$ that has exactly one in-neighbour $j \in \mathcal{N}_k^-$ which is reachable from $\mathcal{U}$. Therefore, we want to find which in-neighbours of $k$ are 
	reachable from $\mathcal{U}$. One of the ways to do this, is to use Dijkstra's single source shortest path (SSSP) algorithm \cite{Dijkstra1959}, \cite{Orlin2010}. This algorithm computes the shortest paths (i.e., paths of minimum length) from a given source node $s$ to every other node in the graph, and returns an `infinite' distance for each node which is not reachable from $s$. If we apply the SSSP algorithm to each node in $\mathcal{U}$, we obtain all nodes in $\mathcal{V}$ that are reachable from $\mathcal{U}$. Dijkstra's SSSP algorithm has time complexity $\mathcal{O}(n + e)$, where $n = \abs{\mathcal{V}}$ and $e = \abs{\mathcal{E}}$ \cite{Orlin2010}, and therefore we can find all nodes reachable from $\mathcal{U}$ in time complexity $\mathcal{O}(un + ue)$, where $u = \abs{\mathcal{U}}$. Once we know the nodes in $\mathcal{V}$ that are reachable from $\mathcal{U}$, we can simply check whether there exists exactly one $j \in \mathcal{N}_k^-$ that is reachable from $\mathcal{U}$. In particular, this shows that the graph simplification process can be implemented in polynomial time since both operations 1 and 2 can be implemented in polynomial time, and the graph simplification process executes at most $n$ operations of type 1 and 2 (if applied in this order).
\end{remark}

\begin{example}
	\label{example5}
	Consider the graph $\mathcal{G} = (\mathcal{V},\mathcal{E})$ in Figure \ref{fig:graph5} and define $\mathcal{U} := \{2\}$ and $\mathcal{W} := \{5,6\}$. The goal of this example is to apply the graph simplification process to obtain a derived graph and derived vertex set. After this simplification, it will be easy to check left-invertibility of $T_{\mathcal{W},\mathcal{U}}(z;G)$. 
	
	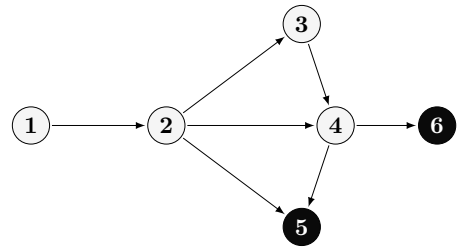
\begin{figure}[h!]
		\centering
		\scalebox{0.9}{
			\begin{tikzpicture}[scale=1]
			\node[VertexStyle1,draw] (1) at (0,0) {$\bf 1$};
			\node[VertexStyle1,draw] (2) at (2,0) {$\bf 2$};
			\node[VertexStyle1,draw] (3) at (4,1.5) {$\bf 3$};
			\node[VertexStyle1,draw] (4) at (4.5,0) {$\bf 4$};
			\node[VertexStyle2,draw] (5) at (4,-1.5) {$\bf 5$};
			\node[VertexStyle2,draw] (6) at (6,0) {$\bf 6$};
			\draw[-latex] (1) -- (2);
			\draw[-latex] (2) -- (4);
			\draw[-latex] (2) -- (3);
			\draw[-latex] (2) -- (5);
			\draw[-latex] (3) -- (4);
			\draw[-latex] (4) -- (5);
			\draw[-latex] (4) -- (6);
			\end{tikzpicture}
		}
		\caption{Graph $\mathcal{G}$ with nodes $\mathcal{W}$ colored black.}
		\label{fig:graph5}
	\end{figure} 
	
	First, note that both nodes $5$ and $6$ do not have outgoing edges, so at the moment we cannot apply operation 1. However, we observe that node $6$ has exactly one in-neighbour (node $4$) that is reachable from $\mathcal{U}$. Consequently, we can replace node $6$ by node $4$ in $\mathcal{W}$ (see Figure \ref{fig:graph6}).
	
	\begin{figure}[h!]
		\centering
		\scalebox{0.9}{
			\begin{tikzpicture}[scale=1]
			\node[VertexStyle1,draw] (1) at (0,0) {$\bf 1$};
			\node[VertexStyle1,draw] (2) at (2,0) {$\bf 2$};
			\node[VertexStyle1,draw] (3) at (4,1.5) {$\bf 3$};
			\node[VertexStyle2,draw] (4) at (4.5,0) {$\bf 4$};
			\node[VertexStyle2,draw] (5) at (4,-1.5) {$\bf 5$};
			\node[VertexStyle1,draw] (6) at (6,0) {$\bf 6$};
			\draw[-latex] (1) -- (2);
			\draw[-latex] (2) -- (4);
			\draw[-latex] (2) -- (3);
			\draw[-latex] (2) -- (5);
			\draw[-latex] (3) -- (4);
			\draw[-latex] (4) -- (5);
			\draw[-latex] (4) -- (6);
			\end{tikzpicture}
		}
		\caption{Graph with nodes $\mathcal{W}$, obtained by applying operation 2 to node $6$.}
		\label{fig:graph6}
	\end{figure}
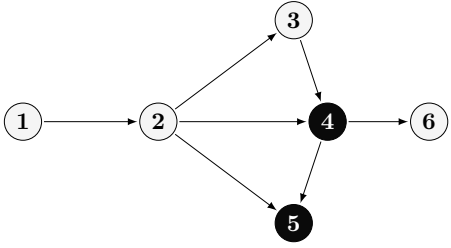
	
	To follow up, we see that node $4$ has outgoing edges, which we can remove by applying operation 1, see Figure \ref{fig:graph7}.
	
	\begin{figure}[h!]
		\centering
		\scalebox{0.9}{
			\begin{tikzpicture}[scale=1]
			\node[VertexStyle1,draw] (1) at (0,0) {$\bf 1$};
			\node[VertexStyle1,draw] (2) at (2,0) {$\bf 2$};
			\node[VertexStyle1,draw] (3) at (4,1.5) {$\bf 3$};
			\node[VertexStyle2,draw] (4) at (4.5,0) {$\bf 4$};
			\node[VertexStyle2,draw] (5) at (4,-1.5) {$\bf 5$};
			\node[VertexStyle1,draw] (6) at (6,0) {$\bf 6$};
			\draw[-latex] (1) -- (2);
			\draw[-latex] (2) -- (4);
			\draw[-latex] (2) -- (3);
			\draw[-latex] (2) -- (5);
			\draw[-latex] (3) -- (4);
			\end{tikzpicture}
		}
		\caption{Graph with nodes $\mathcal{W}$, obtained by applying operation 1 to node $4$.}
		\label{fig:graph7}
	\end{figure}
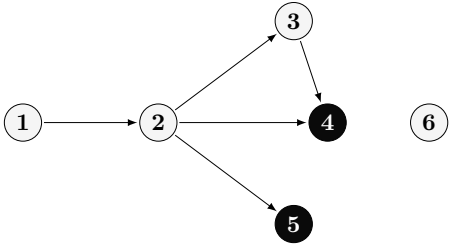
	
	Subsequently, node $5$ has exactly one in-neighbour that is (trivially) reachable from $\mathcal{U}$. Therefore, we replace vertex $5$ by $2$ in $\mathcal{W}$. Next, we can remove all outgoing edges of node $2$ using operation 1. These result of these two operations is depicted in Figure \ref{fig:graph8}.
	
	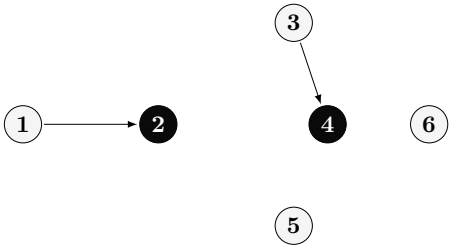
\begin{figure}[h!]
		\centering
		\scalebox{0.9}{
			\begin{tikzpicture}[scale=1]
			\node[VertexStyle1,draw] (1) at (0,0) {$\bf 1$};
			\node[VertexStyle2,draw] (2) at (2,0) {$\bf 2$};
			\node[VertexStyle1,draw] (3) at (4,1.5) {$\bf 3$};
			\node[VertexStyle2,draw] (4) at (4.5,0) {$\bf 4$};
			\node[VertexStyle1,draw] (5) at (4,-1.5) {$\bf 5$};
			\node[VertexStyle1,draw] (6) at (6,0) {$\bf 6$};
			\draw[-latex] (1) -- (2);
			\draw[-latex] (3) -- (4);
			\end{tikzpicture}
		}
		\caption{Derived graph $\mathcal{D}(\mathcal{G})$ with derived vertex set $\mathcal{D}(\mathcal{W})$ (in black), obtained by applying operation 2 to node $5$ and operation 1 to node $2$.}
		\label{fig:graph8}
	\end{figure}
	Note that nodes $2$ and $4$ do not have any outgoing edges. Moreover, the in-neighbour $3$ of node $4$ is not reachable from node $2$, so we cannot use operation 2 to node $4$. In addition, operation 2 cannot be applied to node $2$ since $2 \in \mathcal{U}$. Therefore, the graph simplification process terminates. We conclude that the graph $\mathcal{D}(\mathcal{G})$ in Figure \ref{fig:graph8} is a \emph{derived graph} of $\mathcal{G}$, whereas the vertex set $\mathcal{D}(\mathcal{W}) = \{2,4\}$ is a \emph{derived vertex set} of $\mathcal{W}$ (with respect to $\mathcal{U}$). This example shows the strength of the graph simplification process in the following way: since $\mathcal{U} \subseteq \mathcal{D}(\mathcal{W})$, we conclude by Lemma \ref{lemmanecsys} that $T_{\mathcal{D}(\mathcal{W}),\mathcal{U}}(z;G)$ is left-invertible for all $G(z) \in \mathcal{A}(\mathcal{D}(\mathcal{G}))$. However, by Lemma \ref{lemmaprocess1} and Lemma \ref{lemmaprocess2}, we immediately see that $T_{\mathcal{W},\mathcal{U}}(z;G)$ is left-invertible for all $G(z) \in \mathcal{A}(\mathcal{G})$. This suggests that the graph simplification process is a promising tool to study left-invertibility of transfer matrices (and hence, to study identifiability of dynamical networks).  
\end{example}

To summarize, we have seen that it is possible to remove the outgoing edges of nodes in $\mathcal{W}$ and to `move' the nodes in $\mathcal{W}$ closer to $\mathcal{U}$ if certain conditions are satisfied. Since left-invertibility is preserved by both operations due to Lemmas \ref{lemmaprocess1} and \ref{lemmaprocess2}, we see that left-invertibility of $T_{\mathcal{W},\mathcal{U}}(z;G)$ for all $G(z) \in \mathcal{A}(\mathcal{G})$ is equivalent to the left-invertibility of $T_{\mathcal{D}(\mathcal{W}),\mathcal{U}}(z;G)$ for all $G(z) \in \mathcal{A}(\mathcal{D}(\mathcal{G}))$. Using Lemma \ref{lemmanecsys}, this shows that the condition $\mathcal{U} \subseteq \mathcal{D}(\mathcal{W})$ is \emph{sufficient} for the left-invertibility of $T_{\mathcal{W},\mathcal{U}}(z;G)$. Remarkably, the condition $\mathcal{U} \subseteq \mathcal{D}(\mathcal{W})$ turns out to be also \emph{necessary} for left-invertibility of $T_{\mathcal{W},\mathcal{U}}(z;G)$. This is stated more formally in the following theorem, which is one of the main results of this paper.

\begin{theorem}
	\label{maintheorem}
	Consider a directed graph $\mathcal{G} = (\mathcal{V},\mathcal{E})$ and let $\mathcal{U}, \mathcal{W} \subseteq \mathcal{V}$. Let $\mathcal{D}(\mathcal{W})$ be any derived vertex set of $\mathcal{W}$ with respect to $\mathcal{U}$. Then $\rank T_{\mathcal{W},\mathcal{U}}(z;G) = \abs{\mathcal{U}}$ for all matrices $G(z) \in \mathcal{A}(\mathcal{G})$ if and only if $\mathcal{U} \subseteq \mathcal{D}(\mathcal{W})$.
\end{theorem}

Before we prove Theorem \ref{maintheorem}, we need some auxiliary results. Consider a directed graph $\mathcal{G} = (\mathcal{V},\mathcal{E})$, let $n = \abs{\mathcal{V}}$, $s = \abs{\mathcal{E}}$, and index the edges as $\mathcal{E} = \{e_1,e_2,\dots,e_s\}$. We associate with each edge $e \in \mathcal{E}$ an indeterminate $\mathsf{g}_{e}$. Moreover, we define the s-dimensional vector 
\begin{equation*}
\mathsf{g} := \begin{pmatrix}
\mathsf{g}_{e_1} & \mathsf{g}_{e_2} & \dots & \mathsf{g}_{e_s} \end{pmatrix}^\top,
\end{equation*}
which we call the \emph{indeterminate vector} of $\mathcal{G}$. Next, we define the $n \times n$ matrix $\mathsf{G}$ as
\begin{equation*}
\mathsf{G}_{ji} = \begin{cases*}
\mathsf{g}_{e_k} & \text{if} $e_k = (i,j)$ for some $k$ \\
0 & otherwise.
\end{cases*}
\end{equation*}
We emphasize that not all entries of $\mathsf{G}$ are indeterminates, but some are fixed zeros. Note that we write $\mathsf{G}$ in sans-serif font, to clearly distinguish between $\mathsf{G}$ and a \emph{fixed} rational matrix $G(z)$. It is clear that the determinants of square submatrices of $I - \mathsf{G}$ are real polynomials in the indeterminate entries of $\mathsf{G}$, i.e., in the indeterminate vector $\mathsf{g}$. Hence, the entries of the adjugate of $I- \mathsf{G}$ are real polynomials in $\mathsf{g}$. We state the following basic lemma, which gives conditions under which an entry of $\adj(I- \mathsf{G})$ is a \emph{nonzero} polynomial.

\begin{lemma}
	\label{lemmapolynomial}
	Consider a directed graph $\mc{G} = (\mc{V},\mc{E})$ and let $i,j \in \mathcal{V}$. Let $\mathsf{g}$ and $\mathsf{G}$ be the indeterminate vector and matrix of $\mc{G}$, respectively, and define $\mathsf{A} := \adj(I-\mathsf{G})$. Then $\mathsf{A}_{ji}$ is a nonzero polynomial in $\mathsf{g}$ if and only if there exists a path from $i$ to $j$. 
\end{lemma}
Lemma \ref{lemmapolynomial} follows from Proposition 5.1 of \cite{Hendrickx2018}. Next, we state the following basic result on polynomials. 
\begin{proposition}
	\label{propositionpolynomials}
	Consider $k$ nonzero real polynomials $p_i(x)$, where $i = 1,2,...,k$ and $x = (x_1,x_2,\dots,x_n)$. There exists an $\bar{x} \in \mathbb{R}^n$ such that $p_i(\bar{x}) \neq 0$ for all $i = 1,2,...,k$.   
\end{proposition}
\begin{remark}
	\label{remarkpolynomials}
	Without loss of generality, we can assume that $\bar{x}$ in Proposition \ref{propositionpolynomials} has only nonzero coordinates. Indeed, by continuity, if $p_i(\bar{x}) \neq 0$ for $i = 1,2,...,k$, there exists an open ball $B(\bar{x})$ around $\bar{x}$ in which $p_i(x) \neq 0$ for all $i = 1,2,...,k$ and all $x \in B(\bar{x})$. Clearly, this open ball contains a point with only nonzero coordinates.
\end{remark}
Finally, we require a proposition on rational matrices.
\begin{proposition}
	\label{propositionrationalmatrix}
	Let $A(z)$ be an $m \times n$ rational matrix and assume that each row of $A(z)$ contains at least one nonzero entry. There exists a vector $b \in \mathbb{R}^n$ such that each entry of $A(z)b$ is a nonzero rational function.
\end{proposition}
The proof of Proposition \ref{propositionrationalmatrix} follows simply from induction on the number of rows of $A(z)$ and is therefore omitted. With these results in place, we are ready to prove Theorem \ref{maintheorem}.

\begin{myproof}[Proof of Theorem \ref{maintheorem}]
	Let $\mathcal{D}(\mathcal{G})$ and $\mathcal{D}(\mathcal{W})$ be a derived graph and derived vertex set with respect to $\mathcal{U}$ obtained from the graph simplification process. To prove the `if' statement, suppose that $\mathcal{U} \subseteq \mathcal{D}(\mathcal{W})$. By Corollary \ref{lemmanecsys} we find that $\rank T_{\mathcal{D}(\mathcal{W}),\mathcal{U}}(z;G) = \abs{\mathcal{U}}$ for all $G(z) \in \mathcal{A}(\mathcal{D}(\mathcal{G}))$. By consecutive application of Lemmas \ref{lemmaprocess1} and \ref{lemmaprocess2}, we conclude that $\rank T_{\mathcal{W},\mathcal{U}}(z;G) = \abs{\mathcal{U}}$ for all $G(z) \in \mathcal{A}(\mathcal{G})$. 
	
	Conversely, to prove the `only if' statement, suppose that $\mathcal{U} \not\subseteq \mathcal{D}(\mathcal{W})$. We want to show that 
	\begin{equation*}
	\rank T_{\mathcal{D}(\mathcal{W}),\mathcal{U}}(z;G) < \abs{\mathcal{U}} \text{ for some } G(z) \in \mathcal{A}(\mathcal{D}(\mathcal{G})).
	\end{equation*}
	Since $\mathcal{U} \not\subseteq \mathcal{D}(\mathcal{W})$, the set $\bar{\mathcal{U}} := \mathcal{U} \setminus \mathcal{D}(\mathcal{W})$ is nonempty. Furthermore, as $\mathcal{D}(\mathcal{G})$ and $\mathcal{D}(\mathcal{W})$ result from the graph simplification process, it is clear that nodes in $\mathcal{D}(\mathcal{W})$ do not have outgoing edges. In addition, each node in the set $\bar{\mathcal{W}} := \mathcal{D}(\mathcal{W}) \setminus \mathcal{U}$ has either \emph{zero} or \emph{at least two} in-neighbours that are reachable from $\mathcal{U}$. As nodes in $\mathcal{D}(\mathcal{W}) \cap \mathcal{U}$ have no outgoing edges, this means that each node in $\bar{\mathcal{W}}$ has either zero or at least two in-neighbours that are reachable \emph{from} $\bar{\mathcal{U}}$. Finally, we assume that the nodes in $\mathcal{U}$ do not have any incoming edges, which is without loss of generality by Remark \ref{remarkUincoming}.
	
	The idea of the proof is to show that $T_{\mathcal{D}(\mathcal{W}),\bar{\mathcal{U}}}(z;G) b = 0$, for some to-be-determined network matrix $G(z) \in \mathcal{A}(\mathcal{D}(\mathcal{G}))$ and nonzero vector $b$. Hence, $\rank T_{\mathcal{D}(\mathcal{W}),\bar{\mathcal{U}}}(z;G) < |\bar{\mathcal{U}}|$ and since $T_{\mathcal{D}(\mathcal{W}),\bar{\mathcal{U}}}$ is a submatrix of $T_{\mathcal{D}(\mathcal{W}),\mathcal{U}}$, it will then immediately follow that $\rank T_{\mathcal{D}(\mathcal{W}),\mathcal{U}}(z;G) < \abs{\mathcal{U}}$.
	
	We investigate a row $T_{w,\bar{\mathcal{U}}}(z;G)$ of the transfer matrix $T_{\mathcal{D}(\mathcal{W}),\bar{\mathcal{U}}}(z;G)$ and we distinguish two cases, namely the case that $w \in \mathcal{D}(\mathcal{W}) \cap \mathcal{U}$ and the case that $w \in \bar{\mathcal{W}}$. First, suppose that $w \in \mathcal{D}(\mathcal{W}) \cap \mathcal{U}$. This implies that $w \in \mathcal{U}$. Recall that the nodes in $\mathcal{U}$ do not have any incoming edges. Consequently, there are no paths from $v$ to $w$ for any $v \in \bar{\mathcal{U}}$. We conclude from Lemma 3 of \cite{vandenHof2013} that $T_{wv}(z;G) = 0$ for all $G(z) \in \mathcal{A}(\mathcal{D}(\mathcal{G}))$. Therefore, $T_{w,\bar{\mathcal{U}}}(z;G) = 0$ for all $G(z) \in \mathcal{A}(\mathcal{D}(\mathcal{G}))$. Obviously, this implies that $T_{w,\bar{\mathcal{U}}}(z;G) b = 0$ for all $G(z) \in \mathcal{A}(\mathcal{D}(\mathcal{G}))$ and all real vectors $b$.
	
	Next, we consider the second case in which $w \in \bar{\mathcal{W}}$. Let $\mathsf{G}$ denote the indeterminate matrix of $\mathcal{D}(\mathcal{G})$. In addition, define $\mathsf{A} := \adj(I - \mathsf{G})$. Then, we have
	\begin{subequations}
		\begin{align}
		(I - \mathsf{G}) \mathsf{A} &= \det(I-\mathsf{G})I \\
		(I - \mathsf{G})_{\bar{\mathcal{W}},\mathcal{V}} \mathsf{A}_{\mathcal{V},\bar{\mathcal{U}}} &= 0, \label{varGT}
		\end{align}
	\end{subequations}
	where \eqref{varGT} follows from the fact that $\bar{\mathcal{U}}$ and $\bar{\mathcal{W}}$ are disjoint. Recall that nodes in $\bar{\mathcal{W}}$ do not have any outgoing edges, and therefore $(I - \mathsf{G})_{\bar{\mathcal{W}},\bar{\mathcal{W}}} = I$. This means that we can rewrite \eqref{varGT} as
	\begin{equation}
	\label{varGT2}
	\mathsf{A}_{\bar{\mathcal{W}},\bar{\mathcal{U}}} =
	\mathsf{G}_{\bar{\mathcal{W}},\bar{\mathcal{W}}^c} \mathsf{A}_{\bar{\mathcal{W}}^c,\bar{\mathcal{U}}},
	\end{equation}
	where we recall that $\bar{\mathcal{W}}^c := \mathcal{V} \setminus \bar{\mathcal{W}}$. Note that for $j \in \bar{\mathcal{W}}^c$, the column $\mathsf{G}_{\bar{\mathcal{W}},j}$ is equal to $0$ if $j$ is not an in-neighbour of any node in $\bar{\mathcal{W}}$. In addition, for any $j \in \bar{\mathcal{W}}^c$, the row $\mathsf{A}_{j,\bar{\mathcal{U}}}$ equals $0$ if there is no path from $\bar{\mathcal{U}}$ to $j$ (by Lemma \ref{lemmapolynomial}). Therefore, we can rewrite \eqref{varGT2} as
	\begin{equation}
	\label{relationA}
	\mathsf{A}_{\bar{\mathcal{W}},\bar{\mathcal{U}}} =
	\mathsf{G}_{\bar{\mathcal{W}},\mathcal{N}} \mathsf{A}_{\mathcal{N},\bar{\mathcal{U}}},
	\end{equation}
	where $\mathcal{N} \subseteq \bar{\mathcal{W}}^c$ is characterized by the following property: we have $j \in \mathcal{N}$ if and only if $j$ is an in-neighbour of a node in $\bar{\mathcal{W}}$ and there is a path from $\bar{\mathcal{U}}$ to $j$. By definition of the adjugate, the entries of $\mathsf{A}_{\mathcal{N},\bar{\mathcal{U}}}$ are polynomials in the indeterminate entries of $\mathsf{G}$. We claim that the indeterminate entries of $\mathsf{G}_{\bar{\mathcal{W}},\mathcal{N}}$ do not appear in any entry of $\mathsf{A}_{\mathcal{N},\bar{\mathcal{U}}}$, that is, $\mathsf{A}_{\mathcal{N},\bar{\mathcal{U}}}$ is \emph{independent} of the indeterminate entries of $\mathsf{G}_{\bar{\mathcal{W}},\mathcal{N}}$. For the sake of clarity, we postpone the proof of this claim to the end. For now, we assume that $\mathsf{A}_{\mathcal{N},\bar{\mathcal{U}}}$ is independent of the indeterminate entries of $\mathsf{G}_{\bar{\mathcal{W}},\mathcal{N}}$.
	
By definition, there is a path from $\bar{\mathcal{U}}$ to each node in $\mathcal{N}$. Let $\mathcal{N} = \{n_1,n_2,\dots, n_{r}\}$, where $r = \abs{\mathcal{N}}$. Then, for each node $n_i \in \mathcal{N}$, there exists a node $u_i \in \bar{\mathcal{U}}$ such that $\mathsf{A}_{n_i,u_i}$ is a nonzero polynomial in the indeterminate entries of $\mathsf{G}$ (by Lemma \ref{lemmapolynomial}). We emphasize that $u_i$ and $u_j$ are not necessarily distinct. We focus on the $r$ nonzero polynomials 
	\begin{equation}
	\label{rpolynomials}
	\mathsf{A}_{n_1,u_1}, \mathsf{A}_{n_2,u_2}, \dots, \mathsf{A}_{n_r,u_r}.
	\end{equation}
	The idea is to apply Proposition \ref{propositionpolynomials} and Remark \ref{remarkpolynomials} to these $r$ polynomials. By Remark \ref{remarkpolynomials}, we can substitute nonzero real numbers for the indeterminate entries of $\mathsf{G}$ such that all $r$ polynomials \eqref{rpolynomials} evaluate to nonzero real numbers. Since the polynomials \eqref{rpolynomials} are independent of the indeterminate entries of $\mathsf{G}_{\bar{\mathcal{W}},\mathcal{N}}$, we do not have to fix the entries of $\mathsf{G}_{\bar{\mathcal{W}},\mathcal{N}}$. In addition, it is possible to substitute \emph{strictly proper functions} in $z$ for the indeterminate entries of $\mathsf{G}$ (except for entries of $\mathsf{G}_{\bar{\mathcal{W}},\mathcal{N}}$) such that the polynomials \eqref{rpolynomials} evaluate to \emph{nonzero} rational functions. Indeed, one can simply choose all indeterminate entries of $\mathsf{G}$ as nonzero real numbers as before, and then divide all of these real numbers by $z$. 
	
	To summarize the progress so far, we have substituted strictly proper functions for the indeterminate entries of $\mathsf{G}$ (except for the entries of $\mathsf{G}_{\bar{\mathcal{W}},\mathcal{N}}$) such that the polynomials \eqref{rpolynomials} evaluate to \emph{nonzero} rational functions. Note that this implies that the matrix $\mathsf{A}_{\mathcal{N},\bar{\mathcal{U}}}$ evaluates to a rational matrix, which we denote by $A_{\mathcal{N},\bar{\mathcal{U}}}(z)$ from now on. Since each row of $A_{\mathcal{N},\bar{\mathcal{U}}}(z)$ contains a nonzero rational function, by Proposition \ref{propositionrationalmatrix} there exists a nonzero real vector $b$ such that $A_{\mathcal{N},\bar{\mathcal{U}}}(z) b$ has only \emph{nonzero} rational entries.  
	
	Subsequently, we will choose the indeterminate entries of $\mathsf{G}_{\bar{\mathcal{W}},\mathcal{N}}$ such that $\mathsf{G}_{\bar{\mathcal{W}},\mathcal{N}} A_{\mathcal{N},\bar{\mathcal{U}}}(z) b = 0$. Recall that the nodes in $\bar{\mathcal{W}}$ either have zero or at least two in-neighbours from the set $\mathcal{N}$. If a node $w \in \bar{\mathcal{W}}$ has no in-neighbours, then $\mathsf{G}_{w,\mathcal{\mathcal{N}}} = 0$, and therefore clearly $\mathsf{G}_{w,\mathcal{N}} A_{\mathcal{N},\bar{\mathcal{U}}}(z) b = 0$. If a node $w \in \bar{\mathcal{W}}$ has at least two in-neighbours, say $n_1,n_2,\dots,n_p \in \mathcal{N}$, then we substitute \emph{strictly proper} functions for the indeterminate entries $\mathsf{G}_{w,n_1}, \mathsf{G}_{w,n_2}, \dots, \mathsf{G}_{w,n_p}$ so that $\mathsf{G}_{w,\mathcal{N}} A_{\mathcal{N},\bar{\mathcal{U}}}(z) b = 0$. Note that this is possible since the vector $A_{\mathcal{N},\bar{\mathcal{U}}}(z) b$ has only \emph{nonzero} rational entries. To conclude, we have substituted strictly proper functions for the indeterminate entries of $\mathsf{G}$ which yields a matrix which we denote by $G(z)$. The adjugate of $I-G(z)$ is denoted by $A(z) = \adj(I-G(z))$. We have shown that $G_{\bar{\mathcal{W}},\mathcal{N}}(z) A_{\mathcal{N},\bar{\mathcal{U}}}(z) b = 0$. By \eqref{relationA}, this yields $A_{\bar{\mathcal{W}},\bar{\mathcal{U}}}(z) b = 0$. Note that $\det (I - G(z))$ is nonzero since all nonzero entries of $G$ are strictly proper functions. Therefore, 
	\begin{equation*}
	T(z;G) = \frac{1}{\det(I-G(z))} A(z),
	\end{equation*}
	from which we find that $T_{\bar{\mathcal{W}},\bar{\mathcal{U}}}(z;G) b = 0$. Consequently, $T_{\mathcal{D}(\mathcal{W}),\bar{\mathcal{U}}}(z;G) b = 0$, and $\rank T_{\mathcal{D}(\mathcal{W}),\bar{\mathcal{U}}}(z;G) < \abs{\bar{\mathcal{U}}}$. Therefore, we conclude that $\rank T_{\mathcal{D}(\mathcal{W}),\mathcal{U}}(z;G) < \abs{\mathcal{U}}$. We still have to show that $G(z)$ is admissible, i.e., $G(z) \in \mathcal{A}(\mathcal{D}(\mathcal{G}))$. Since the indeterminate matrix $\mathsf{G}$ is consistent with the graph $\mathcal{D}(\mathcal{G})$ and we substituted (nonzero) strictly proper functions for each indeterminate entry of $\mathsf{G}$, the matrix $G(z)$ readily satisfies Properties P1 and P2. In addition, since all nonzero entries of $G(z)$ are strictly proper, we obtain $
	\lim_{z \to \infty} I - G(z) = I$, and hence, $G(z)$ also satisfies Property P3. We conclude that $\rank T_{\mathcal{D}(\mathcal{W}),\mathcal{U}}(z;G) < \abs{\mathcal{U}}$ for some $G(z) \in \mathcal{A}(\mathcal{D}(\mathcal{G}))$. Finally, by consecutive application of Lemmas \ref{lemmaprocess1} and \ref{lemmaprocess2}, we conclude that $\rank T_{\mathcal{W},\mathcal{U}}(z;G) < \abs{\mathcal{U}}$ for some $G(z) \in \mathcal{A}(\mathcal{G})$. 
	
	Recall that we have so far assumed that $\mathsf{A}_{\mathcal{N},\bar{\mathcal{U}}}$ is \emph{independent} of the indeterminate entries of $\mathsf{G}_{\bar{\mathcal{W}},\mathcal{N}}$. It remains to be shown that this is true. To this end, label the nodes in $\mathcal{V}$ such that $\mathsf{G}$ can be written as
	\begin{subequations}
	\begin{align}
	\mathsf{G} &= \begin{pmatrix}
	\mathsf{G}_{\bar{\mathcal{W}}^c,\bar{\mathcal{W}}^c} & \mathsf{G}_{\bar{\mathcal{W}}^c,\bar{\mathcal{W}}} \\
	\mathsf{G}_{\bar{\mathcal{W}},\bar{\mathcal{W}}^c} & \mathsf{G}_{\bar{\mathcal{W}},\bar{\mathcal{W}}}
	\end{pmatrix} \\
	&= \begin{pmatrix}
	\mathsf{G}_{\bar{\mathcal{W}}^c,\bar{\mathcal{W}}^c} & 0 \\
	\mathsf{G}_{\bar{\mathcal{W}},\bar{\mathcal{W}}^c} & 0
	\end{pmatrix}, \label{eqWc2}
	\end{align}
	\end{subequations}
	where \eqref{eqWc2} follows from the fact that nodes in $\bar{\mathcal{W}}$ have no outgoing edges. This implies that 
	\begin{equation*}
	I - \mathsf{G} = \begin{pmatrix}
	I - \mathsf{G}_{\bar{\mathcal{W}}^c,\bar{\mathcal{W}}^c} & 0 \\
	-\mathsf{G}_{\bar{\mathcal{W}},\bar{\mathcal{W}}^c} & I
	\end{pmatrix},
	\end{equation*}
	and therefore
	\begin{equation}
	\label{formulaadj}
	\mathsf{A} = \adj(I - \mathsf{G}) = \begin{pmatrix}
	\adj(I - \mathsf{G}_{\bar{\mathcal{W}}^c,\bar{\mathcal{W}}^c}) & 0 \\
	* & *
	\end{pmatrix}.
	\end{equation}
	Since the entries of $\mathsf{G}_{\bar{\mathcal{W}}^c,\bar{\mathcal{W}}^c}$ are independent of the indeterminate entries of $\mathsf{G}_{\bar{\mathcal{W}},\bar{\mathcal{W}}^c}$, we conclude from \eqref{formulaadj} that the matrix $\mathsf{A}_{\bar{\mathcal{W}}^c,\bar{\mathcal{W}}^c} = \adj(I - \mathsf{G}_{\bar{\mathcal{W}}^c,\bar{\mathcal{W}}^c})$ is independent of the indeterminate entries of $\mathsf{G}_{\bar{\mathcal{W}},\bar{\mathcal{W}}^c}$. Now, to prove the claim, note that $\bar{\mathcal{U}}$ and $\bar{\mathcal{W}}$ are disjoint by definition, and therefore $\bar{\mathcal{U}} \subseteq \bar{\mathcal{W}}^c$. In addition, we have $\mathcal{N} \subseteq \bar{\mathcal{W}}^c$. Therefore, the matrix $\mathsf{A}_{\mathcal{N},\bar{\mathcal{U}}}$ is a \emph{submatrix} of $\mathsf{A}_{\bar{\mathcal{W}}^c,\bar{\mathcal{W}}^c}$. Furthermore, we see that $\mathsf{G}_{\bar{\mathcal{W}},\mathcal{N}}$ is a submatrix of $\mathsf{G}_{\bar{\mathcal{W}},\bar{\mathcal{W}}^c}$ by using the fact that $\mathcal{N} \subseteq \bar{\mathcal{W}}^c$. We conclude that the entries of $\mathsf{A}_{\mathcal{N},\bar{\mathcal{U}}}$ are independent of the indeterminate entries of $\mathsf{G}_{\bar{\mathcal{W}},\mathcal{N}}$.
\end{myproof}

\section{Identifiability and graph simplification}
\label{sectionidentifiabilitysimplification}

In this section we use Theorem \ref{maintheorem} to provide solutions to the identifiability problems introduced in Section \ref{sectionproblem}. Specifically, the following theorem follows from Theorem \ref{maintheorem} and Lemma \ref{lemma1} and states necessary and sufficient graph-theoretic conditions for identifiability of $(i,\mathcal{N}^+_i)$. 
\begin{theorem}
	\label{mainresultiNi}
	Consider a directed graph $\mathcal{G} = (\mathcal{V},\mathcal{E})$, let $i \in \mathcal{V}$ and $\mathcal{C} \subseteq \mathcal{V}$. Moreover, let $\mathcal{D}(\mathcal{C})$ be any derived vertex set of $\mathcal{C}$ with respect to $\mathcal{N}^+_i$. Then $(i,\mathcal{N}^+_i)$ is identifiable from $\mathcal{C}$ in $\mathcal{G}$ if and only if $\mathcal{N}^+_i \subseteq \mathcal{D}(\mathcal{C})$.
\end{theorem}

\begin{example}
	Consider the graph in Figure \ref{fig:graph5}. We wonder whether $(1,\mathcal{N}^+_1)$ is identifiable. Note that we have $\mathcal{N}_1^+ = \{2\}$. The set of measured nodes is $\mathcal{C} = \{5,6\}$. As shown in Example \ref{example5}, a derived vertex set of $\mathcal{C}$ with respect to $\mathcal{N}_1^+$ is given by $\mathcal{D}(\mathcal{C}) = \{2,4\}$. Since $\{2\} \subseteq \mathcal{D}(\mathcal{C})$, we conclude by Theorem \ref{mainresultiNi} that $(1,\mathcal{N}_1^+)$ is identifiable. In other words, we can uniquely reconstruct $G_{21}(z)$ from the transfer matrix $CT(z;G)$. This approach shows the strength of our approach. Indeed, note that to check identifiability, we do not have to verify Definition \ref{def1} directly. Also, we do not have to compute $CT(z;G) = C(I-G(z))^{-1}$ and verify its rank for all $G(z) \in \mathcal{A}(\mathcal{G})$, which is required to check the condition of Lemma \ref{lemma1}.
\end{example}

By definition of the graph simplification process, we have that $|\mathcal{D}(\mathcal{C})| \leq |\mathcal{C}|$. Hence, it follows from Theorem \ref{mainresultiNi} that identifiability of $(i,\mathcal{N}_i^+)$ implies that the number of measured nodes is greater or equal to the number of out-neighbours of node $i$. 
\begin{corollary}
	\label{cormainNi}
	Consider a directed graph $\mathcal{G} = (\mathcal{V},\mathcal{E})$, let $i \in \mathcal{V}$ and $\mathcal{C} \subseteq \mathcal{V}$. If $(i,\mathcal{N}^+_i)$ is identifiable from $\mathcal{C}$ in $\mathcal{G}$ then $|\mathcal{N}^+_i| \leq |\mathcal{C}|$.
\end{corollary}

The next result gives necessary and sufficient graph-theoretic conditions under which the entire graph $\mathcal{G}$ is identifiable. This result is a corollary of Theorem \ref{mainresultiNi} but is stated as a theorem due to its importance.

\begin{theorem}
	\label{mainresultG}
	Consider a directed graph $\mathcal{G} = (\mathcal{V},\mathcal{E})$ and let $\mathcal{C} \subseteq \mathcal{V}$. Then $\mathcal{G}$ is identifiable from $\mathcal{C}$ if and only if for all $i \in \mathcal{V}$, we have $\mathcal{N}^+_i \subseteq \mathcal{D}(\mathcal{C})$, where $\mathcal{D}(\mathcal{C})$ is any derived vertex set of $\mathcal{C}$ with respect to $\mathcal{N}^+_i$. 
\end{theorem} 

We emphasize that the derived set $\mathcal{D}(\mathcal{C})$ of $\mathcal{C}$ depends on the choice of neighbour set $\mathcal{N}_i^+$, and hence, for each node $i \in \mathcal{V}$ we have to compute the derived set of $\mathcal{C}$ with respect to $\mathcal{N}_i^+$.

\section{Comparison to results based on constrained vertex-disjoint paths}
\label{sectionresultscvdp}

In the previous section we established necessary and sufficient graph-theoretic conditions for the identifiability of respectively $(i,\mathcal{N}_i^+)$ and $\mathcal{G}$. The purpose of the current section is to compare these results to the ones based on so-called \emph{constrained vertex-disjoint paths} \cite{vanWaarde2018}. We first recall the definition in what follows. 

\begin{definition}
	Let $\mathcal{G} = (\mathcal{V},\mathcal{E})$ be a directed graph. Consider a set of $m$ vertex-disjoint paths in $\mathcal{G}$ with starting nodes $\bar{\mathcal{U}} \subseteq \mathcal{V}$ and end nodes $\bar{\mathcal{W}} \subseteq \mathcal{V}$. We say that the set of vertex-disjoint paths is \emph{constrained} if it is the \emph{only} set of $m$ vertex-disjoint paths from $\bar{\mathcal{U}}$ to $\bar{\mathcal{W}}$.
\end{definition} 

Next, let $\mathcal{U}, \mathcal{W} \subseteq \mathcal{V}$ be disjoint subsets of vertices. We say that there exists a constrained set of $m$ vertex-disjoint paths \emph{from} $\mathcal{U}$ \emph{to} $\mathcal{W}$ if there exists a constrained set of $m$ vertex-disjoint paths in $\mathcal{G}$ with starting nodes $\bar{\mathcal{U}} \subseteq \mathcal{U}$ and end nodes $\bar{\mathcal{W}} \subseteq \mathcal{W}$. In the case that $\mathcal{U} \cap \mathcal{W} \neq \emptyset$, we say that there is a constrained set of $m$ vertex-disjoint paths from $\mathcal{U}$ to $\mathcal{W}$ if there exists a constrained set of $\max\{0,m - \abs{\mathcal{U} \cap \mathcal{W}}\}$ vertex-disjoint paths from $\mathcal{U} \setminus \mathcal{W}$ to $\mathcal{W} \setminus \mathcal{U}$.

\begin{remark}
	Note that for a set of $m$ vertex-disjoint paths from $\mathcal{U}$ to $\mathcal{W}$ to be constrained, we do not require the existence of a unique set of $m$ vertex-disjoint paths from $\mathcal{U}$ to $\mathcal{W}$. In fact, we only require the existence of a unique set of vertex-disjoint paths between the \emph{starting nodes} $\bar{\mathcal{U}}$ of the paths and the \emph{end nodes} $\bar{\mathcal{W}}$. We will illustrate the definition of constrained vertex-disjoint paths in Example \ref{example2}.
\end{remark}

\begin{remark}
	The notion of constrained vertex-disjoint paths is strongly related to the notion of \emph{constrained matchings} in bipartite graphs \cite{Hershkowitz1993}. In fact, a constrained matching can be seen as a special case of a constrained set of vertex-disjoint paths where all paths are of length one. 
\end{remark}

\begin{example}
	\label{example2}
	Consider the graph $\mathcal{G} = (\mathcal{V},\mathcal{E})$ in Figure \ref{fig:graph2}. Moreover, consider the subsets of vertices $\mathcal{U} := \{2,3\}$ and $\mathcal{W} := \{6,7,8\}$.
	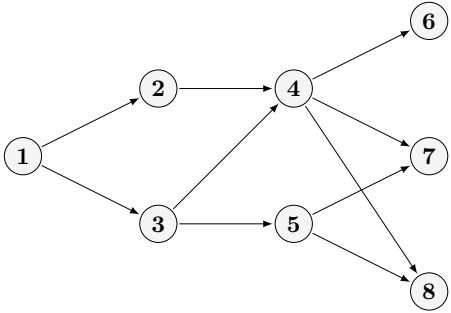
\begin{figure}[h!]
		\centering
		\scalebox{0.9}{
			\begin{tikzpicture}[scale=1]
			\node[VertexStyle1,draw] (1) at (0,0) {$\bf 1$};
			\node[VertexStyle1,draw] (2) at (2,1) {$\bf 2$};
			\node[VertexStyle1,draw] (3) at (2,-1) {$\bf 3$};
			\node[VertexStyle1,draw] (4) at (4,1) {$\bf 4$};
			\node[VertexStyle1,draw] (5) at (4,-1) {$\bf 5$};
			\node[VertexStyle1,draw] (6) at (6,2) {$\bf 6$};
			\node[VertexStyle1,draw] (7) at (6,0) {$\bf 7$};
			\node[VertexStyle1,draw] (8) at (6,-2) {$\bf 8$};
			\draw[-latex] (1) -- (2);
			\draw[-latex] (1) -- (3);
			\draw[-latex] (2) -- (4);
			\draw[-latex] (3) -- (4);
			\draw[-latex] (3) -- (5);
			\draw[-latex] (4) -- (6);
			\draw[-latex] (4) -- (7);
			\draw[-latex] (4) -- (8);
			\draw[-latex] (5) -- (7);
			\draw[-latex] (5) -- (8);
			\end{tikzpicture}
		}
		\caption{Graph used in Example \ref{example2}.}
		\label{fig:graph2}
	\end{figure} 
	Clearly, the paths $\{(2,4), (4,6)\}$ and $\{(3,5), (5,7)\}$ form a set of two vertex-disjoint paths from $\mathcal{U}$ to $\mathcal{W}$. In fact, this set of vertex-disjoint paths is \emph{constrained} since there does not exist another set of two vertex-disjoint paths from $\bar{\mathcal{U}} = \{2,3\}$ to $\bar{\mathcal{W}} = \{6,7\}$. Therefore, there exists a constrained set of two vertex-disjoint paths from $\mathcal{U}$ to $\mathcal{W}$. Note that there are also other sets of vertex-disjoint paths from $\mathcal{U}$ to $\mathcal{W}$. For example, the paths $\{(2,4),(4,7)\}$ and $\{(3,5),(5,8)\}$ also form a set of two vertex-disjoint paths. However, this set of vertex-disjoint paths is \emph{not} constrained. To see this, note that we have another set of vertex-disjoint paths from $\bar{\mathcal{U}} = \{2,3\}$ to $\bar{\mathcal{W}} = \{7,8\}$, namely the set consisting of the paths $\{(2,4),(4,8)\}$ and $\{(3,5),(5,7)\}$. 
\end{example}

In the following theorem, we recall the main result presented in \cite{vanWaarde2018}, which relates the notion of constrained vertex-disjoint paths and identifiability of $(i,\mathcal{N}_i^+)$.

\begin{theorem}
	\label{theoremnecsys}
	Consider a directed graph $\mathcal{G} = (\mathcal{V},\mathcal{E})$, let $i \in \mathcal{V}$ and $\mathcal{C} \subseteq \mathcal{V}$. If there exists a constrained set of $| \mathcal{N}_i^+ |$ vertex-disjoint paths from $\mathcal{N}_i^+$ to $\mathcal{C}$ then $(i,\mathcal{N}_i^+)$ is identifiable from $\mathcal{C}$.
\end{theorem}

The proof of Theorem \ref{theoremnecsys} can be found in \cite{vanWaarde2018} (see Theorem 13). A natural question to ask is whether the condition given in Theorem \ref{theoremnecsys} is also \emph{necessary} for identifiability. It turns out that this is not the case, as demonstrated next.

\begin{example}
	\label{example4}
	In this example, we revisit the graph $\mathcal{G} = (\mathcal{V},\mathcal{E})$ in Figure \ref{fig:graph5}. Suppose that we are interested in the identifiability of $(1,\mathcal{N}_1^+)$, i.e., in the identifiability of the transfer function corresponding to the edge $(1,2)$. The set of measured nodes is given by $\mathcal{C} = \{5,6\}$. The purpose of this example is to show that Theorem \ref{theoremnecsys} is not necessary, i.e., we have to show that $(1,\mathcal{N}_1^+)$ is identifiable even though there does not exist a constrained set of one (vertex-disjoint) path from $\mathcal{N}_1^+$ to $\mathcal{C}$.
	
	Note that $\mathcal{N}_1^+ = \{2\}$ and that there are three different paths from $2$ to $5$. In addition, there are two different paths from node $2$ to node $6$. This implies that there does not exist a constrained set of one (vertex-disjoint) path from $\mathcal{N}_1^+$ to $\mathcal{C}$. Nonetheless, we can show that $(1,\mathcal{N}_1^+)$ is identifiable. The easiest way to show this is by noting that we already proved in Example \ref{example5} that $\mathcal{N}_1^+ \subseteq \mathcal{D}(\mathcal{C})$, where $\mathcal{D}(\mathcal{C})$ is a derived vertex set of $\mathcal{C}$. Hence, by Theorem \ref{mainresultiNi} we conclude that $(1,\mathcal{N}_1^+)$ is identifiable. However, to gain a bit more insight we will prove identifiability of $(1,\mathcal{N}_1^+)$ by inspection of the transfer matrix $T_{\mathcal{C},\mathcal{N}_1^+}(z;G)$. For any $G(z) \in \mathcal{A}(\mathcal{G})$, we obtain
	\begin{equation}
	\label{TCN}
	T_{\mathcal{C},\mathcal{N}_1^+} = \begin{pmatrix}
	G_{52} + G_{54}(G_{42} + G_{43}G_{32}) \\
	G_{64}(G_{42} + G_{43}G_{32})
	\end{pmatrix},
	\end{equation}
	where we omitted the argument $z$. If $G_{42} + G_{43}G_{32} \neq 0$ then $G_{64}(G_{42} + G_{43}G_{32}) \neq 0$ and therefore $\rank T_{\mathcal{C},\mathcal{N}_1^+} = 1$. If $G_{42} + G_{43}G_{32} = 0$, we see that $G_{52} + G_{54}(G_{42} + G_{43}G_{32}) = G_{52} \neq 0$ so also in this case $\rank T_{\mathcal{C},\mathcal{N}_1^+} = 1$. We conclude that   $\rank T_{\mathcal{C},\mathcal{N}_1^+} = 1$ for all admissible network matrices, which means that $(1,\mathcal{N}_1^+)$ is identifiable by Lemma \ref{lemma1}. 
\end{example}

Example \ref{example4} also gives some intuition for the fact that Theorem \ref{theoremnecsys} is not necessary for identifiability. Indeed, the condition based on constrained vertex-disjoint paths guarantees that a \emph{square} submatrix of $T_{\mathcal{C},\mathcal{N}_i^+}(z;G)$ is invertible \emph{for all} admissible $G$, where the columns and rows of this submatrix are indexed by the starting nodes and end nodes of the paths, respectively \cite{vanWaarde2018}. However, as can be seen from \eqref{TCN}, the matrix $T_{\mathcal{C},\mathcal{N}_i^+}(z;G)$ might be left-invertible for all admissible $G$, even though there \emph{does not exist} a square $\abs{\mathcal{N}_i^+} \times \abs{\mathcal{N}_i^+}$ submatrix of $T_{\mathcal{C},\mathcal{N}_i^+}(z;G)$ that is invertible for all admissible $G$. 
In general, the particular square submatrix of $T_{\mathcal{C},\mathcal{N}_i^+}(z;G)$ that is invertible \emph{depends} on the network matrix $G$. Interestingly, we can use the general theory developed in this paper to show that the condition of Theorem \ref{theoremnecsys} is necessary and sufficient in the special case that $T_{\mathcal{C},\mathcal{N}_i^+}(z;G)$ is \emph{square itself} (a proof is given in Appendix \ref{appendixtheoremnecsuf}).

\begin{theorem}
	\label{theoremnecsuf}
	Consider a directed graph $\mathcal{G} = (\mathcal{V},\mathcal{E})$. Let $i \in \mathcal{V}$ and $\mathcal{C} \subseteq \mathcal{V}$ be such that $\abs{\mathcal{C}} = \abs{\mathcal{N}_i^+}$. Then, $(i,\mathcal{N}_i^+)$ is identifiable if and only if there exists a constrained set of $\abs{\mathcal{N}_i^+}$ vertex-disjoint paths from $\mathcal{N}_i^+$ to $\mathcal{C}$.
\end{theorem}

The main message of this section is that the conditions in terms of constrained vertex-disjoint paths \cite{vanWaarde2018} are only necessary and sufficient in the special case that $\abs{\mathcal{N}_i^+} = \abs{\mathcal{C}}$. This case is quite particular, especially if one is interested in identifiability of the entire network. In the latter situation, Theorem \ref{theoremnecsuf} can only be applied if the number of out-neighbours of \emph{each node} is equal to the number of measured nodes, which is very restrictive. Therefore, we conclude that the necessary and sufficient conditions for identifiability based on graph simplification are much more general. Additional advantages of the conditions based on the graph simplification process are that they are conceptually simpler and appealing from computational point of view, see Remark \ref{remarkcomputation}. 

\section{Conclusions}
\label{sectionconclusions}

In this paper we have considered the problem of identifiabi\-lity of dynamical networks for which interactions between nodes are modelled by transfer functions. We have been interested in graph-theoretic conditions for two identifiability problems. First, we wanted to find conditions under which the transfer functions of all outgoing edges of a given node are identifiable. Secondly, we have been interested in conditions under which all transfer functions in the network are identifiable. It is known that these problems are equivalent to the left-invertibility of certain transfer matrices \emph{for all} networked matrices associated with the graph \cite{Hendrickx2018}, \cite{vanWaarde2018}. However, the downside of such rank conditions is that it is not clear how to \emph{check} the rank of a transfer matrix for an \emph{infinite} number of network matrices.

Therefore, as our first contribution, we have provided a necessary and sufficient graph-theoretic condition under which a transfer matrix has full column rank \emph{for all} network matrices. To this end, we have introduced a new concept called the \emph{graph simplification process}. The idea of this process is to apply simplifying operations to the graph, after which left-invertibility can be verified by simply checking a set inclusion. Based on the graph simplification process, we have given necessary and sufficient conditions for identifiability. Notably, we have shown that our conditions can be verified by polynomial time algorithms. Finally, we have shown that our results generalize existing sufficient conditions based on constrained vertex-disjoint paths \cite{vanWaarde2018}.

It is interesting to observe that our topological conditions for global identifiability are quite different from the path-based conditions for \emph{generic} identifiability \cite{Hendrickx2018}. This is analogous to the \emph{controllability} literature, where it was shown that weak structural controllability can be characterized in terms of maximal matchings \cite{Liu2011}, while strong structural controllability was characterized using a (different) graph-theoretic concept called zero forcing \cite{Monshizadeh2014}.  

For future work, it would be interesting to consider a \emph{minimum sensor placement} problem. The goal in such a problem is to find sets of measured nodes of minimum cardinality such that the entire network is identifiable. 

\appendix

\subsection{Proof of Lemma \ref{lemmanecsys}}
\label{appendixlemmanecsys}

\begin{myproof}[Proof of Lemma \ref{lemmanecsys}]
	Without loss of generality, we assume that the nodes in $\mathcal{W}$ do not have outgoing edges (see Lemma \ref{lemmaprocess1}). Since $\mathcal{U} \subseteq \mathcal{W}$, the nodes in $\mathcal{U}$ do not have outgoing edges. We now relabel the nodes in $\mathcal{G}$ such that $G(z) \in \mathcal{A}(\mathcal{G})$ can be written as
	\begin{equation*}
	G = \begin{pmatrix}
	G_{\mathcal{U},\mathcal{U}} & G_{\mathcal{U},\mathcal{U}^c} \\
	G_{\mathcal{U}^c,\mathcal{U}} & G_{\mathcal{U}^c,\mathcal{U}^c}
	\end{pmatrix} = \begin{pmatrix}
	0 & G_{\mathcal{U},\mathcal{U}^c} \\
	0 & G_{\mathcal{U}^c,\mathcal{U}^c}
	\end{pmatrix},
	\end{equation*}
	where we omitted the argument $z$, and where the zeros are present due to the fact that nodes in $\mathcal{U}$ do not have outgoing edges. Consequently, we obtain 
	\begin{equation*}
	T = (I - G)\inv = \begin{pmatrix}
	I & -G_{\mathcal{U},\mathcal{U}^c} \\
	0 & I-G_{\mathcal{U}^c,\mathcal{U}^c}
	\end{pmatrix}\inv = \begin{pmatrix}
	I & * \\
	0 & *
	\end{pmatrix},
	\end{equation*}
	and therefore, $T_{\mathcal{U},\mathcal{U}} = I$. Hence, $T_{\mathcal{U},\mathcal{U}}$ has full rank for all $G(z) \in \mathcal{A}(\mathcal{G})$ and we conclude that $T_{\mathcal{W},\mathcal{U}}$ has rank $\abs{\mathcal{U}}$ for all $G(z) \in \mathcal{A}(\mathcal{G})$.
\end{myproof}

\subsection{Proof of Theorem \ref{theoremnecsuf}}
\label{appendixtheoremnecsuf}

To prove Theorem \ref{theoremnecsuf}, we will first state two lemmas. Under the assumption that $\abs{\mathcal{U}} = \abs{\mathcal{W}}$, the following lemma asserts that the existence of a set of constrained vertex-disjoint paths from $\mathcal{U}$ to $\mathcal{W}$ is preserved by operation 1.

\begin{lemma}
\label{lemmacp1}
Let $\mathcal{G} = (\mathcal{V},\mathcal{E})$ be a directed graph and consider $\mathcal{U}, \mathcal{W} \subseteq \mathcal{V}$ such that $\abs{\mathcal{U}} = \abs{\mathcal{W}}$. Moreover, let $\bar{\mathcal{G}} = (\mathcal{V},\bar{\mathcal{E}})$ be the graph obtained from $\mathcal{G}$ by removing all outgoing edges of the nodes in $\mathcal{W}$. There exists a constrained set of $\abs{\mathcal{U}}$ vertex-disjoint paths from $\mathcal{U}$ to $\mathcal{W}$ in $\mathcal{G}$ if and only if there exists a constrained set of $\abs{\mathcal{U}}$ vertex-disjoint paths from $\mathcal{U}$ to $\mathcal{W}$ in $\bar{\mathcal{G}}$.
\end{lemma}

\begin{myproof}
The lemma follows from the following important observation: if $\abs{\mathcal{U}} = \abs{\mathcal{W}}$, then a set of $\abs{\mathcal{U}}$ vertex-disjoint paths from $\mathcal{U}$ to $\mathcal{W}$ \emph{does not contain} any outgoing edge of a node in $\mathcal{W}$. Indeed, if a path $\mathcal{P}$ from $\mathcal{U}$ to $\mathcal{W}$ in such a set of vertex-disjoint paths contains an edge $(w,v)$, where $w \in \mathcal{W}$ and $v \in \mathcal{V}$, then the path $\mathcal{P}$ contains at least two vertices in $\mathcal{W}$ (namely $w$ and the end node). This means that $\mathcal{P}$ is contained in a set of at most $\abs{\mathcal{U}} - 1$ vertex disjoint paths from $\mathcal{U}$ to $\mathcal{W}$. However, this is a contradiction since we assumed that $\mathcal{P}$ was contained in a set of $\abs{\mathcal{U}}$ vertex-disjoint paths from $\mathcal{U}$ to $\mathcal{W}$.

Next, we prove the `if' statement. Suppose that there exists a constrained set $\mathcal{S}$ of $\abs{\mathcal{U}}$ vertex-disjoint paths from $\mathcal{U}$ to $\mathcal{W}$ in $\bar{\mathcal{G}}$. Then $\mathcal{S}$ is also a set of $\abs{\mathcal{U}}$ vertex-disjoint paths from $\mathcal{U}$ to $\mathcal{W}$ in $\mathcal{G}$. We want to prove that $\mathcal{S}$ is constrained (in the graph $\mathcal{G}$). Therefore, suppose on the contrary that there exists another set $\bar{\mathcal{S}}$ of $\abs{\mathcal{U}}$ vertex-disjoint paths from $\mathcal{U}$ to $\mathcal{W}$ in $\mathcal{G}$. By the above discussion, we know that no path in $\bar{\mathcal{S}}$ contains an outgoing edge of a node in $\mathcal{W}$. Therefore, $\bar{\mathcal{S}}$ is a set of $\abs{\mathcal{U}}$ vertex-disjoint paths from $\mathcal{U}$ to $\mathcal{W}$ in $\bar{\mathcal{G}}$. As such, we conclude that $\bar{\mathcal{S}} = \mathcal{S}$. In other words, $\mathcal{S}$ is a constrained set of $\abs{\mathcal{U}}$ vertex-disjoint paths from $\mathcal{U}$ to $\mathcal{W}$ in $\mathcal{G}$.

Conversely, to prove the `only if' statement, suppose that there exists a constrained set $\mathcal{S}$ of $\abs{\mathcal{U}}$ vertex-disjoint paths from $\mathcal{U}$ to $\mathcal{W}$ in $\mathcal{G}$. Again, by the previous discussion we know that no path in $\mathcal{S}$ contains an outgoing edge of a node in $\mathcal{W}$. Therefore, $\mathcal{S}$ is also a constrained set of $\abs{\mathcal{U}}$ vertex-disjoint paths from $\mathcal{U}$ to $\mathcal{W}$ in $\bar{\mathcal{G}}$. This proves the lemma.
\end{myproof}

The following lemma relates the existence of a constrained set of vertex-disjoint paths and the \emph{second} graph operation.

\begin{lemma}
\label{lemmacp2}
Consider a directed graph $\mathcal{G} = (\mathcal{V},\mathcal{E})$ and let $\mathcal{U},\mathcal{W} \subseteq \mathcal{V}$. Suppose that $k \in \mathcal{W} \setminus \mathcal{U}$ has exactly one in-neighbour $j \in \mathcal{N}_k^-$ that is reachable from $\mathcal{U}$. Then there exists a constrained set of $\abs{\mathcal{U}}$ vertex-disjoint paths from $\mathcal{U}$ to $\mathcal{W}$ if and only if there exists a constrained set of $\abs{\mathcal{U}}$ vertex-disjoint paths from $\mathcal{U}$ to $\bar{\mathcal{W}}$ with $\bar{\mathcal{W}} := (\mathcal{W} \setminus \{k\}) \cup \{j\}$.
\end{lemma}

\begin{myproof}
We will first show that $\mathcal{S}$ is a set of $\abs{\mathcal{U}}$ vertex-disjoint paths from $\mathcal{U}$ to $\mathcal{W}$ if and only if $\bar{\mathcal{S}}$ is a set of $\abs{\mathcal{U}}$ vertex-disjoint paths from $\mathcal{U}$ to $\bar{\mathcal{W}}$, where $\bar{\mathcal{S}}$ will be specified. 

Suppose that $\mathcal{S}$ is a set of $\abs{\mathcal{U}}$ vertex disjoint paths from $\mathcal{U}$ to $\mathcal{W}$. Consider the path $\mathcal{P} \in \mathcal{S}$ that goes from $\mathcal{U}$ to $k$. Since $j \in \mathcal{N}_k^-$ is the only in-neighbour of $k$ that is reachable from $\mathcal{U}$, we obtain $(j,k) \in \mathcal{P}$. This means that $\bar{\mathcal{P}} := \mathcal{P} \setminus (j,k)$ is a path from $\mathcal{U}$ to $j$. Clearly, $\bar{\mathcal{S}} := (\mathcal{S} \setminus \mathcal{P}) \cup \bar{\mathcal{P}}$ is a set of $\abs{\mathcal{U}}$ vertex-disjoint paths from $\mathcal{U}$ to $\bar{\mathcal{W}}$.

Conversely, suppose that $\bar{\mathcal{S}}$ is a set of $\abs{\mathcal{U}}$ vertex-disjoint paths from $\mathcal{U}$ to $\bar{\mathcal{W}}$. Consider the path $\bar{\mathcal{P}} \in \bar{\mathcal{S}}$ that goes from $\mathcal{U}$ to $j \in \bar{\mathcal{W}}$. Since $j \in \mathcal{N}_k^-$ is the only in-neighbour of $k$ that is reachable from $\mathcal{U}$, the path $\bar{\mathcal{P}}$ does not pass through the vertex $k$. Consequently, $\mathcal{P} := \bar{\mathcal{P}} \cup (j,k)$ is a path from $\mathcal{U}$ to $k$. Again using the fact that $j$ is the only in-neighbour of $k$ that is reachable from $\mathcal{U}$, we see that no path in $\bar{\mathcal{S}}$ passes through the vertex $k$. This implies that $\mathcal{S} := (\bar{\mathcal{S}} \setminus \bar{\mathcal{P}}) \cup \mathcal{P}$ is a set of $\abs{\mathcal{U}}$ vertex-disjoint paths from $\mathcal{U}$ to $\mathcal{W}$.

To conclude, we have shown that $\mathcal{S}$ is a set of $\abs{\mathcal{U}}$ vertex-disjoint paths from $\mathcal{U}$ to $\mathcal{W}$ if and only if $\bar{\mathcal{S}}$ is a set of $\abs{\mathcal{U}}$ vertex-disjoint paths from $\mathcal{U}$ to $\bar{\mathcal{W}}$, where the set $\bar{\mathcal{S}}$ is defined as $\bar{\mathcal{S}}:= (\mathcal{S} \setminus \mathcal{P}) \cup \bar{\mathcal{P}}$. This implies that $\mathcal{S}$ is a \emph{constrained} set of $\abs{\mathcal{U}}$ vertex-disjoint paths from $\mathcal{U}$ to $\mathcal{W}$ if and only if $\bar{\mathcal{S}}$ is a \emph{constrained} set of $\abs{\mathcal{U}}$ vertex-disjoint paths from $\mathcal{U}$ to $\bar{\mathcal{W}}$.
\end{myproof}


\begin{myproof}[Proof of Theorem \ref{theoremnecsuf}]
The `if' statement follows from Theorem \ref{theoremnecsys}. To prove the `only if' part, suppose that $(i,\mathcal{N}_i^+)$ is identifiable. By Theorem \ref{mainresultiNi}, $\mathcal{N}_i^+ \subseteq \mathcal{D}(\mathcal{C})$, where $\mathcal{D}(\mathcal{C})$ is a derived vertex set of $\mathcal{C}$ with respect to $\mathcal{N}_i^+$. In fact, we obtain $\mathcal{N}_i^+ = \mathcal{D}(\mathcal{C})$ as $\abs{\mathcal{N}_i^+} = \abs{\mathcal{C}}$. Let $\mathcal{D}(\mathcal{G})$ denote the associated derived graph of $\mathcal{G}$. By definition, there exists a constrained set of $\abs{\mathcal{N}_i^+}$ vertex-disjoint paths from $\mathcal{N}_i^+$ to $\mathcal{D}(\mathcal{C})$ in $\mathcal{D}(\mathcal{G})$ (see Section \ref{sectionresultscvdp}). By consecutive application of Lemmas \ref{lemmacp1} and \ref{lemmacp2}, we conclude that there exists a constrained set of $\abs{\mathcal{N}_i^+}$ vertex-disjoint paths from $\mathcal{N}_i^+$ to $\mathcal{C}$ in the graph $\mathcal{G}$.
\end{myproof}

\bibliographystyle{IEEEtran}
\bibliography{MyRef}

\vfill

\end{document}